\RequirePackage{etoolbox}
\csdef{input@path}{{style/}{graphics/}}
\documentclass[aos,MSNbibl,seceqn,nameyear,dvips]{arximspdf}
\usepackage{graphicx}
\doi{10.1214/14-AOS1288}
\volume{43}
\issue{2}
\pubyear{2015}
\firstpage{620}
\lastpage{651}
\docsubty{FLA}

\makeatletter

\def\cal{\mathcal}

\def\min{\mathrm{min}}
\def\mes{\mathrm{mes}}
\def\R{\mathbb{R}}
\def\cI{{\cal I}}
\def\cJ{{\cal J}}

\newcommand{\hf}{{\hat f}}
\newcommand{\hcF}{{\hat{\cal F}}}
\newcommand{\tf}{{\tilde f}}
\newcommand{\hth}{{\hat\theta}}
\newcommand{\bg}{\mathbf{g}}
\newcommand{\bA}{\mathbf{A}}
\newcommand{\bT}{\mathbf{T}}
\newcommand{\bU}{\mathbf{U}}
\newcommand{\bb}{\mathbf{b}}
\newcommand{\bc}{\mathbf{c}}
\newcommand{\bff}{\mathbf{f}}
\newcommand{\bz}{\mathbf{z}}
\newcommand{\bd}{\mathbf{d}}
\newcommand{\mJ}{{{m}_J}}
\newcommand{\cS}{{\cal S}}
\newcommand{\cF}{{\cal F}}
\newcommand{\cG}{{\cal G}}
\newcommand{\hcG}{{\hat{\cal G}}}
\newcommand{\thetam}{{\bolds\theta}}
\newcommand{\deltam}{{\bolds\delta}}
\newcommand{\betam}{{\bolds\beta}}
\newcommand{\etam}{{\bolds\eta}}
\newcommand{\tmu}{{\tilde\mu}}
\newcommand{\mum}{{\bolds\mu}}
\newcommand{\babff}{{\bar\bff}}
\newcommand{\hbff}{{\hat\bff}}
\newcommand{\bathetam}{{\bar\thetam}}
\newcommand{\hthetam}{{\hat\thetam}}
\newcommand{\argmin}{\mathop{\arg\min}}

\newcommand{\bzero}{{\mathbf 0}}
\newcommand{\bone}{{\mathbf 1}}
\newcommand{\ve}{\varepsilon}

\newproclaim{definition}{Definition}
\newtheorem{lemma}{Lemma}

\newtheorem{theorem}{Theorem}
\newproclaim{remark}{Remark}

\newproclaim{Example}{Example}

\makeatother

\begin{document}
\begin{frontmatter}

\title{Asymptotics for in-sample density forecasting}
\runtitle{In-sample density forecasting}

\begin{aug}
\author[A]{\fnms{Young K.}~\snm{Lee}\thanksref{T1}\ead[label=e1]{youngklee@kangwon.ac.kr}}, %
\author[B]{\fnms{Enno}~\snm{Mammen}\thanksref{T2}\ead[label=e2]{mammen@math.uni-heidelberg.de}},
\author[C]{\fnms{Jens P.}~\snm{Nielsen}\thanksref{T3}\ead[label=e3]{jens.nielsen.1@city.ac.uk}} %
\and\break 
\author[D]{\fnms{Byeong U.}~\snm{Park}\corref{}\ead[label=e4]{bupark@stats.snu.ac.kr}\ead[label=e5]{bupark2000@gmail.com}\thanksref{T4}}
\runauthor{Lee, Mammen, Nielsen and Park}
\affiliation{Kangwon National University, Universit\"at Heidelberg
and Higher School of Economics, Moscow,
Cass Business School, City University London and
Seoul National University}

\thankstext{T1}{Supported by the NRF grant funded by the Korea
government (MEST) (No. 2010-0021396).}
\thankstext{T2}{Supported by the Collaborative Research Center SFB 884
``Political Economy of Reforms,'' financed by the German Science
Foundation (DFG).}
\thankstext{T3}{Supported by
the Institute and Faculty of Actuaries, London, UK.}
\thankstext{T4}{Supported by the NRF grant funded by the Korea
government (MEST)
(No. 2010-0017437).}
\address[A]{Y.~K. Lee\\
Department of Statistics\\
Kangwon National University\\
Chuncheon 200-701\\
Korea\\
\printead{e1}}
\address[B]{E. Mammen\\
Institute fuer Angewandte Mathematik\\
Universitaet Heidelberg\\
Im Neuenheimer Feld 294\\
69120 Heidelberg\\
Germany\\
\printead{e2}}
\address[C]{J. P. Nielsen\\
Cass Business School\\
City University London\\
106 Bunhill Row\\
London EC1Y 8TZ\\
United Kingdom\\
\printead{e3}}
\address[D]{B.~U. Park\\
Department of Statistics\\
Seoul National University\\
Seoul 151-747\\
Korea\\
\printead{e4}\\
\phantom{E-mail:\ }\printead*{e5}\hspace*{48pt}}
\end{aug}

%
\received{\smonth{1} \syear{2014}}
%
\revised{\smonth{9} \syear{2014}}

%
\begin{abstract}
This paper generalizes recent proposals of density forecasting models
and it develops theory for this class of models. In density
forecasting, the density of observations is estimated in regions where
the density is not observed. Identification of the density in such
regions is guaranteed by structural assumptions on the density that
allows exact extrapolation. In this paper, the structural assumption is
made that the density is a product of one-dimensional functions. The
theory is quite general in assuming the shape of the region where the
density is observed. Such models naturally arise when the time point of
an observation can be written as the sum of two terms (e.g., onset and
incubation period of a disease). The developed theory also allows for a
multiplicative factor of seasonal effects. Seasonal effects are present
in many actuarial, biostatistical, econometric and statistical studies.
Smoothing estimators are proposed that are based on backfitting. Full
asymptotic theory is derived for them. A practical example from the
insurance business is given producing a within year budget of reported
insurance claims.
A small sample study supports the theoretical results.
\end{abstract}

%
\end{frontmatter}

\section{Introduction}\label{sec1}

In-sample density forecasting is in this paper defined as forecasting a
structured density in regions where the density is not observed.
This is possible when the density is structured in such a way that all
entering components are estimable in-sample.
Let us, for example, assume that we have one covariate $X$ representing
the start of something; it could be onset of
some infection, underwriting of an insurance contract or the reporting
of an insurance claim, birth of a new member
of a cohort or an employee losing his job in the labour market. Let
then $Y$ represent the development or delay
to some event from this starting point. It could be incubation period
of some disease, development of an insurance
claim, age of a cohort member or time spend looking for a new job. Then
$X+Y$ is the calendar time of the relevant
event. This event is observed if and only if it has already happened
until a calendar time, say $t_0$.
The forecasting exercise is about predicting the density of future
events in calendar times after $t_0$.

The most typical example of a structured density is a simple
multiplicative form studied
by Mammen, Mart{\'\i}nez-Miranda and Nielsen (\citeyear{MamMarNie13}).
The multiplicative density model assumes that $X$ and $Y$ are
independent with smooth densities $f$ and~$g$.
When $f$ and $g$ are estimated by histograms, our in-sample forecasting
approach could be formulated via
a parametric model. This version of in-sample density forecasting is
omnipresent in academic studies as well as
in business forecasting; see Mart{\'\i}nez-Miranda et al. (\citeyear{Maretal13}) for more details and references
in insurance and in statistics of cohort models. Extensions of such
parametric histogram type of models can often be understood
as structured density models modelled via histograms. A structured
density is defined as a known function of
lower-dimensional unknown underlying functions; see \citet{MamNie03} for a formal definition of
generalised structured models. Under the assumption that the model is
true, our forecasts do not extrapolate any parameters or time series
into the future. We therefore call our methodology ``in-sample density
forecasting'': a~structured density estimator forecasting the future
without further assumptions or approximate extrapolations.

Our model is related to deconvolution, but there are two major
differences. First, in our model
one observes not only $X+Y$ but also the summands $X$ and $Y$.
Second, $X$ and $Y$ are only
observed if their sum lies in a certain set, for example, in an
interval $(0,t_0]$. This makes
$X$ and $Y$ be dependent and the estimation problem be an inverse problem. We
will see below that the
first difference leads to rates of convergence that coincide with rates
for the estimation of one-dimensional
functions in the classical nonparametric regression and density
settings. The reason is that our model
consists in a well-posed inverse problem. In contrast, deconvolution is
an ill-posed inverse problem
and allows only poorer rates of convergence.

This paper adds three new contributions to the literature on in-sample
density forecasting. First of all, we define smoothing estimators based
on backfitting and we develop a complete asymptotic distribution theory
for these estimators. Second, we allow for a general class of regions
for which the density is observed. The leading example is a triangle. A
triangle arises in the above examples where the sum of two covariates
is bounded by calendar time. The theoretical discussion in Mammen,
Mart{\'\i}nez-Miranda and Nielsen (\citeyear{MamMarNie13}) was restricted to this case.
But there exist many other important support sets; see, for example,
\citet{KuaNieNie08} for a detailed discussion. Third,
we generalize the forecasting model by modelling a seasonal component.
This is done by introducing an additional multiplicative seasonal
factor into the model. Then we have three one-dimensional density
functions that enter the model and that can be estimated in sample.
Seasonal effects are omnipresent: onset of some disease could be more
likely in the winter than in the summer; new jobs might be less likely
during the summer or they may depend on the business cycle; more auto
insurance claims are reported during the winter, but they might be
bigger on average in the summer; cold winters or hot summers affect
mortality. When a study is running over a few years only and one or two
of those years are not fully observed, data might be too sparse to
leave these two years out of the study. Leaving them in might however
generate bias. The inclusion of seasonality in this paper solves this
type of problems and allow us in general to do well when years are not
fully observed. An illustration producing a within-year budget of
insurance claims is given in the application section.

Classical actuarial methodology does not include seasonal effects.
Budgets are normally carried out manually by highly paid actuaries. The
automatic adjustment of seasonal effects offered by this paper is
therefore potentially cost saving. Insurance companies currently use
the classical chain ladder technique when forecasting future claims.
Classical chain ladder has recently been identified as being the above
mentioned multiplicative histogram in-sample forecasting approach; see
Mart{\'\i}nez-Miranda et al. (\citeyear{Maretal13}). The seasonal
adjustment suggested in this paper is therefore directly implementable
to working routines and processes used by today's nonlife insurance companies.

Recent updates of classical chain ladder include \citet{KuaNieNie09},
\citet{VerNieJes10},
Mart{\'\i}nez-Miranda et al. (\citeyear{Maretal11}) and Mart{\'\i
}nez-Miranda, Nielsen and Verrall (\citeyear{MarNieVer12}). These papers re-\break interpreted
classical chain ladder in modern mathematical statistical terms. The
generalised structured nonparametric model of this paper is a
multiplicative density with three effects. The third seasonal effect is
a function of the covariates of the first two effects. Estimation is
carried out by projecting an unstructured local linear density
estimator, \citet{Nie99}, down on the structure of interest. The
seasonal addition to the multiplicative density model of Mammen, Mart{\'
\i}nez-Miranda and Nielsen (\citeyear{MamMarNie13}) is still a generalised additive
structure, a simple special case of generalised structured models.
Generalised structured models have historically been more studied in
regression than in density estimation. Future developments of our
in-sample density approach will therefore naturally be related to
fundamental regression models; see \citet{LinNie95}, \citet{NieLin98}, \citet{OpsRup97}, \citet{MamLinNie99}, \citet{JiaFanFan10},
Mammen and Park (\citeyear{MamPar05,MamPar06}), \citet{NieSpe05}, \citet{MamNie03},
\citet{YuParMam08}, Lee, Mammen and Park
(\citeyear{LeeMamPar10,LeeMamPar12,LeeMamPar14}),
\citet{ZhaParWan13}, among others.

The paper is structured as follows. Section~\ref{sec2} describes our structured in-sample
density forecasting model, and show that the model is identifiable
(estimable) under weak conditions.
Section~\ref{methodology} explains a new approach to the estimation of the model. Here,
it is assumed
that the data are observed in continuous time and nonparametric
smoothing methods are applied. Section~\ref{sec4}
contains the theoretical properties of our method and Section~\ref{sec5}
considers numerical examples and discusses
the performance of the new approach. The \hyperref[app]{Appendix} contains technical details.

\section{The model}\label{sec2}

We observe a random sample $\{(X_i,Y_i)\dvtx 1\leq i \leq n\}$ from a
density $f$ supported on a subset $\cI$ of a rectangle
$[0,1]^2$.
The density $f(x,y)$ of $(X_i,Y_i)$ is a multiplicative function of
three univariate components, where
the first two are a function of the coordinate $x$ and $y$,
respectively, and the third is a function of the sum of the
two coordinates, $x+y$, and is periodic. Specifically, we consider the following
multiplicative model:
%
\begin{equation}
\label{model} f(x,y) = f_1(x) f_2(y) f_3
\bigl(\mJ(x+y)\bigr),\qquad (x,y) \in\cI,
\end{equation}
where $\mJ(t) =J\mbox{mod}_J(t)$, $\mathrm{ mod}_J(t)= t \mbox{ modulo }
1/J$ for some $J>0$, that is,
$m_J(t)=J(t-l/J)$ for $l/J \le t < (l+1)/J$, $j=0,1,2,\ldots.$ Here,
$f_j$ are unknown
nonnegative functions supported and bounded away from zero on their supports.
We note that $m_J(t)$ always takes values
in $[0,1)$ as $t$ varies on $\R^+$, and that the third component
$f_3(\mJ(\cdot))$
is a periodic function with period $J^{-1}$.

We will prove the identifiability of the functions $f_1$, $f_2$ and
$f_3$ under the constraints that
$\int_0^1 f_1(x) \,dx = \int_0^1 f_2(y) \,dy =1$.
We will do this for two scenarios. In the first case, we assume that
$f_1$, $f_2$ and $f_3$ are smooth functions. Then identification
follows by a simple argument. Our second result does not make use of
smoothness conditions of the component functions. It only requires
conditions on the shape of the set $\cI$. The second result is
important for an understanding of our estimation procedure that is
based on a projection onto the model
(\ref{model}) without using a smoothing procedure for the component functions.

Our first identifiability result makes use of the following conditions:
\begin{longlist}[(A1)]
\item[(A1)] The projections of the set $\cI$
onto the $x$- and $y$-axis equal $[0,1]$.

\item[(A2)] For every $z \in[0,1)$ there exists $(x,y)$ in the
interior of $\cI$ with $m_J(x+y)=z$. Furthermore, for every $x, y \in
(0,1)$ there exist $x'$ and $y'$ with $(x,y')$ and $(x',y)$ in the
interior of $\cI$.
\item[(A3)] The functions $f_{1},f_{2},f_{3}$ are bounded away from
zero and infinity on their supports.
\item[(A4)] The functions $f_{1}$ and $f_{2}$ are differentiable on
$[0,1]$. The function $f_{3}$ is twice differentiable on $[0,1)$.
\item[(A5)] There exist sequences $x_0=0 < x_1 < \cdots< x_k=1$ and $y_0
=1 > y_1 > \cdots > y_k =0$ with $(x,y_j) \in\cI$ for $x_j \leq x \leq x_{j+1}$.
\end{longlist}

\begin{theorem} \label{theoident1} Assume that model (\ref{model})
holds with \textup{(A1)--(A5)}.
Then the functions $f_{1},f_{2},f_{3}$ are identifiable.
\end{theorem}

\begin{remark}
Let $T= \max\{x+y\dvtx (x,y) \in\cI\}$. We note that
the functions $f_j$ are not identifiable in case $J < 1/T$. To see
this, we take
$f_1(u) = f_2(u)= c_1 e^u, f_3(u)=e^u$ with the constant $c_1>0$ chosen
for $f_1=f_2$ to satisfy the constraint $\int_0^1 f_j(u) \,du = 1$.
Consider also $g_1(u)=g_2(u)=c_2 e^{(J+1)u}, g_3(u)=c_1^2/c_2^2$
with the constants $c_2>0$ chosen
for $g_1=g_2$ to satisfy the constraint \mbox{$\int_0^1 g_j(u) \,du = 1$}. In
case $J < 1/T$, we have $m_J(x+y)=J(x+y)$ for all
$(x,y) \in\cI$. This implies that $(f_1,f_2,f_3)$ and $(g_1,g_2,g_3)$
give the same multiplicative density.
In fact, if $J <1/T$, then the assumption (A2) is not fulfilled.
\end{remark}

We now come to our second identifiability result that does not require
smoothness conditions for the functions $f_1$, $f_2$ and $f_3$.
This makes use of the following conditions on the shape of the support
set $\cI$. To introduce conditions on the support set
$\cI$, we let $I_1(y)=\{x\dvtx (x,y) \in\cI\}$, $I_2(x)=\{y\dvtx (x,y) \in
\cI
\}$ and
$I_{3l}(z) = \{x\in[0,1]\dvtx (x,(z+l)/J-x) \in\cI\}$. Below, we assume
that these sets change smoothly as $y, x$ and $z$,
respectively, move. Here, $A \triangle B$ denotes the symmetric
difference of two sets $A$ and $B$ in $\R$, and
$\mes(A)$ the Lebesgue measure of a set $A \subset\R$. Recall the definition
$T= \max\{x+y\dvtx (x,y) \in\cI\}$, and with this define $L(J)$ be the
largest integer that is less than or equal to $TJ$.
\begin{longlist}[(A6)]
\item[(A6)] For $j \in\{1,2,3\}$ there exist partitions $0=a_0^j < \cdots
< a_{L_j}^j=1$ of $[0,1]$ and a function $\kappa\dvtx  [0,1] \to\mathbb
{R}^+$ with $\kappa(x) \to0$ for $x \to0$ such that
(i) for all $u_1, u_2 \in(a_{l-1}^j,a_l^j)$, $\mes[I_j(u_1) \triangle
I_j(u_2)] \le\kappa(|u_1-u_2|), l= 1,\ldots,L_j; j=1,2$;
(ii) for all $u_1, u_2 \in(a_{l-1}^3,a_l^3)$, $\sum_{k=0}^{L(J)}\mes
[I_{3k}(u_1) \triangle I_{3k}(u_2)] \le\kappa(|u_1-u_2|), l= 1,\ldots,L_3$.
Furthermore, it holds that $\mes(I_2(x)) > 0$, $\mes(I_1(y)) > 0$ and
$\sum_{l=0}^{L(J)} \mes[I_{3l}(z)] > 0$ for $x,y\in(0,1)$ and for $z
\in[0,1)$.
\end{longlist}

Assumption (A6) will be used to prove the continuity of some relevant
functions that appear in the technical arguments.
The continuity of a function $\gamma$ implies that $\gamma(x) =0$ for
all $x$ if it is zero almost all $x$.
The assumption allows a finite number of \textit{jumps} in $I_j(u)$ for
$j=1,2$ and $I_{3k}(u)$ as $u$ moves.
For example, suppose that $\cI=\{(x,y)\dvtx 0 \le x \le1, 0 \le y \le
1, x+y \le5/4\}$ and $J=2$. In this case,
$L(J)=2$, and for $k=0,1$ we have $I_{3k}(z) = [0,(z+k)/2]$ for all $z
\in[0,1)$, so that $I_{3k}$ changes smoothly
as $z$ varies on $[0,1)$. However, for $k=2$ we get that
$I_{3k}(z)=[z/2,1]$ for $z \in[0,1/2]$ and $I_{3k}(z)$ is empty for $z
\in(1/2,1)$, thus it changes drastically at $z=1/2$.
In fact, $\lim_{h \rightarrow0}\sum_{k=0}^{L(J)} \mes[I_{3k}(z+h)
\triangle I_{3k}(z-h)] \neq0$ for $z=1/2$.
We note that in this case assumption (A6) holds if we split $[0,1)$
into two partitions, $[0,1/2)$ and $(1/2,1)$.

\begin{figure}

\includegraphics{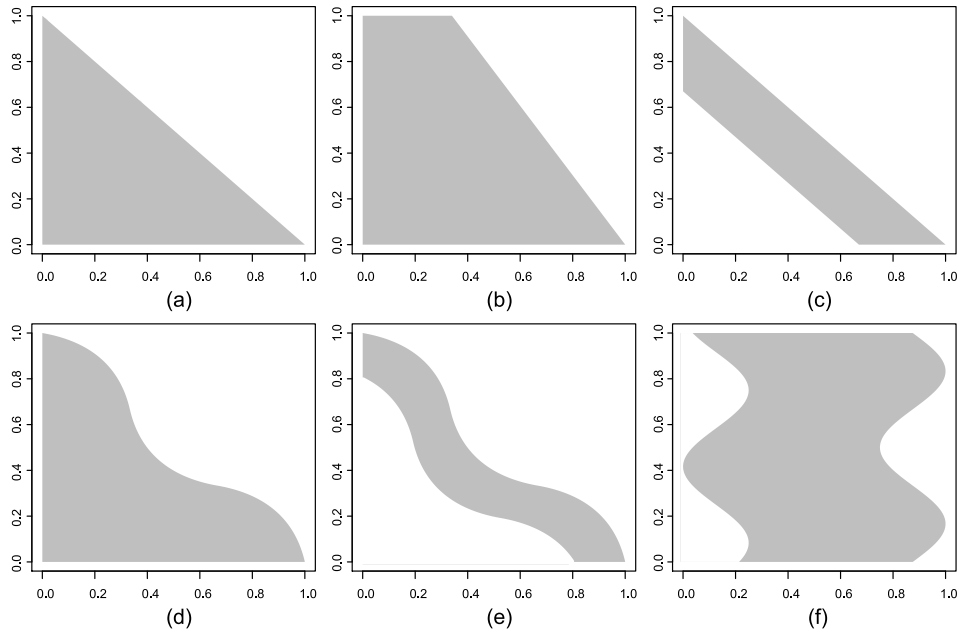}

\caption{Shapes of possible support sets. The horizontal axis
indicates the onset (x) and the vertical the development (y).}
\label{fig:suppset}
\end{figure}

The assumptions (A1), (A2), (A5) and (A6)
accommodate a variety of sets $\cI$ that arise in real applications.
Figure~\ref{fig:suppset} depicts some realistic examples of the set $\cI$ that satisfy
the assumptions. In particular, those sets
of the type in the panels (c) and (e) satisfy (A2) and (A6) if
the maximal vertical or horizontal thickness of the stripe is larger than
the period $1/J$ of the third component function $f_3(m_J(\cdot))$.
In the interpretation of the examples in Figure~\ref{fig:suppset},
we follow the equivalent discussion from Keiding (\citeyear{Ke90}) and Kuang, Nielsen and
Nielsen (\citeyear{KuaNieNie08}). The triangle in Figure~\ref{fig:suppset}(a) is typical for insurance or
mortality when none of the underwriting years or cohorts are fully
run-off. The standard actuarial term ``fully run-off'' means that all
events from that underwriting year or cohort have been observed. In
almost all practical cases of estimating outstanding liabilities,
actuaries stick to the triangle format leaving out fully run-off
underwriting years. While the triangle also appears in mortality
studies, it is common here to leave the fully run-off cohorts in the
study resulting in the support shape given in Figure~\ref{fig:suppset}(b). The support in
Figure~\ref{fig:suppset}(c) arises when the data analyst only considers observations from
the most recent calendar years. While this approach is omnipresent in
practical actuarial science, there is no formal theory or mathematical
models behind these procedures in the actuarial literature. This paper
is therefore an important step toward formalising mathematically
actuarial practise while at the same time improving it. The support
given in Figure~\ref{fig:suppset}(d) and (e) arises when there is a known time
transformation such that time is running at another pace for different
underwriting years or cohort years. While this type of time
transformations are well known in mortality studies are often coined as
versions of accelerated failure time models. Time transformations are
also well known in actuarial science coined as operational time.
However, the academic literature of actuarial science is still
struggling to find a formal definition of what operational time is.
This paper offers one potential solution to this outstanding and
important issue. The last Figure~\ref{fig:suppset}(f) is included to give an impression
of the generality of support structures one could deal with inside our
model approach. Data is missing in the beginning and end of the delay
period, but the model is still valid and in-sample forecasts can be constructed.

The model (\ref{model}) has taken structured density forecasting into a
new territory by leaving the simple multiplicative model. If $f_3$
above was constant (and therefore not in the model) then our model
reduces to the simple multiplicative model analysed in Mart{\'\i
}nez-Miranda et al. (\citeyear{Maretal13}) and Mammen, Mart{\'
\i
}nez-Miranda and Nielsen (\citeyear{MamMarNie13}). These two papers point out that the
simple multiplicative density forecasting model is a continuous version
of a widely used parametric approach corresponding to a structured
histogram version of in-sample density forecasting based on the simple
multiplicative model. The in-sample density forecasting model under
investigation in this paper generalizes the simple multiplicative
approach in an intuitive and simple way including seasonal effects.

In the following theorem, we show that, if there are two multiplicative
representations of the joint density $f$
that agree on almost all points in $\cI$, then the component functions
also agree on almost all points in $[0,1]$.
We will use this result later in the asymptotic analysis of our
estimation procedure.

\begin{theorem} \label{theoident2} Assume that model (\ref{model})
holds with \textup{(A1)--(A3)}, \textup{(A5)}, \textup{(A6)}.
Suppose that $(g_{1},g_{2},g_{3}) $ is a tuple of functions that are
bounded away from zero and infinity with
$\int_0^1 g_1(x) \,dx = \int_0^1 g_2(y) \,dy =1$. Let $\mu_j = \log f_j
- \log g_j$. Assume
that $\mu_1(x) + \mu_2(y) + \mu_3(m_J(x+y))=0$ a.e. on $\cI$.
Then $\mu_j\equiv0$ a.e. on $[0,1]$.
\end{theorem}

\section{Methodology}\label{methodology}

We describe the estimation method for the model (\ref{model}).
We first note that the marginal densities of $X, Y$ and $m_J(X+Y)$
may be zero even if we
assume that the joint density is bounded away from zero. For example,
the marginal densities of $X$ and $Y$ at the
point $u=1$ are zero for the support set $\cI$ given in Figure~\ref{fig:suppset}(a). We
estimate the multiplicative density model on a region where
we observe sufficient data. This means that we exclude the points
$(1,0)$ and $(0,1)$ in the estimation in the case of Figure~\ref{fig:suppset}(a),
and the point $(1,0)$ in the case of Figure~\ref{fig:suppset}(b).
Formally, for a set $S \subset\cI$, let $J_1$ and $J_2$ denote
versions of $I_1$ and $I_2$, respectively,
defined by
$J_1(y)=\{x\dvtx (x,y) \in S\}$ and $J_2(x)=\{y\dvtx (x,y) \in S\}$, and define
$J_{3l}(z)=\{x\dvtx (x,(z+l)/J-x) \in S\}$.
We take an arbitrarily small number $\delta>0$, and find the largest
set $S$ such that
\begin{eqnarray*}
\mes\bigl(J_2(x)\bigr)&\ge&\delta,\qquad \mes\bigl(J_1(y)\bigr)
\ge\delta,\\
\sum_{l=0}^{L(J)}\mes
\bigl(J_{3l}\bigl(m_J(x+y)\bigr)\bigr)&\ge&\delta\qquad \mbox{for
all } (x,y) \in S,
\end{eqnarray*}
where $\mes(A)$ for a set $A$ denotes its length. Such a set is given by
$S=\{(x,y)\dvtx 0 \le x \le1-\delta, 0\le y \le1-\delta, x+y \le1\}$
in the case of Figure~\ref{fig:suppset}(a), and
$S=\{(x,y) \in\cI\dvtx  0 \le x \le1-\delta\}$ in the case of Figure~\ref{fig:suppset}(b),
for example.

We estimate $f_j$ on $S$. Let $S_1$ and $S_2$ be the projections of $S$
onto $x$- and $y$-axis, that is,
$S_1=\{x \in[0,1]\dvtx (x,y) \in S \mbox{ for some } y \in[0,1]\}$,
$S_2=\{y \in[0,1]\dvtx (x,y) \in S \mbox{ for some } x \in[0,1]\}$, and
$S_3=\{m_J(x+y)\dvtx (x,y) \in S\}$. In the case of Figure~\ref{fig:suppset}(a),
$S_1=S_2=[0,1-\delta]$, $S_3=[0,1)$,
but in the case of Figure~\ref{fig:suppset}(b), $S_1=[0,1-\delta]$, $S_2=[0,1]$,
$S_3=[0,1)$. We put the following constraints
on $f_j$:
\[
\int_{S_1}f_1(x) \,dx=\int_{S_2}f_2(y)=1.
\]
This is only for convenience.\vspace*{-1pt}
Now, we define $f_{w,1}(x)=\int_{J_2(x)}f(x,y) \,dy,\break   f_{w,2}(y)=\int_{J_1(y)}f(x,y) \,dx$ and\vspace*{1pt}
$f_{w,3}(z)=\sum_{l=0}^{L(J)} \int_{J_{3l}(z)}f(x,(z+l)/J-x) \,dx$.
Then the model (\ref{model}) gives the following integral equations:
%
\begin{eqnarray}
\label{inteq} %
f_{w,1}(x) &=& f_1(x)\int
_{J_2(x)} f_2(y)f_3\bigl(\mJ(x+y)\bigr)
\,dy,\qquad x \in S_1,\nonumber
\\
f_{w,2}(y) &=& f_2(y)\int_{J_1(y)}
f_1(x)f_3\bigl(\mJ(x+y)\bigr) \,dx, \qquad y \in S_2,
\\
f_{w,3}(z) &=& f_3(z) \sum_{l=0}^{L(J)}
\int_{J_{3l}(z)} f_1(x)f_2\bigl((z+l)/J-x
\bigr) \,dx,\qquad z \in S_3.
\nonumber
\end{eqnarray}
We note that the marginal functions on the left-hand sides of the above
equations are bounded away
from zero on $S_j$. Specifically, $\inf_{u\in S_j} f_{w,j}(u) \ge
\delta\inf_{(x,y)\in\cI}f(x,y)>0$ so that
$f_j$ in the equations are well-defined.

Suppose that we are given a preliminary estimator of the joint density $f$.
Call it $\hf$. We estimate $f_{w,j}$ by $\hf_{w,j}$ that are defined as
$f_{w,j}$, respectively,
with $f$ being replaced by the preliminary estimator $\hf$. Our
proposed estimators of $f_j$, for $j=1,2,3$, are obtained by
replacing $f_{w,j}$ in the integral equations (\ref{inteq}) by $\hf
_{w,j}$, respectively, and solving the resulting
equations for the multiplicative components. Let $\vartheta=\int_S
f(x,y) \,dx \,dy$ and $\hat\vartheta$ be its estimator defined
by $\hat\vartheta=n^{-1}\sum_{i=1}^n I[(X_i,Y_i) \in S]$. Putting the
constraints
%
\begin{eqnarray}
\label{constr} \int_{S_1}\hf_1(x) \,dx &=& \int
_{S_2}\hf_2(y) \,dy =1,
\nonumber
\\[-8pt]
\\[-8pt]
\nonumber
 \int_S
\hf_1(x)\hf_2(y)\hf_3\bigl(m_J(x+y)
\bigr) \,dx \,dy&=&\hat\vartheta,
\end{eqnarray}
they are given as the solution of the following backfitting equations:
%
\begin{eqnarray}
\label{mi} %
\hf_1(x) &=& \hth_1 \cdot
\frac{\hf_{w,1}(x)}{\int_{J_2(x)} \hf
_2(y) \hf
_3(\mJ(x+y)) \,dy},\nonumber
\\
\hf_2(y) &=& \hth_2 \cdot\frac{\hf_{w,2}(y)}{\int_{J_1(y)} \hf
_1(x) \hf
_3(\mJ(x+y)) \,dx},
\\
\hf_3(z) &=& \hth_3 \cdot\frac{\hf_{w,3}(z)}
{\sum_{l=0}^{L(J)}\int_{J_{3l}(z)}\hf_1(x)\hf_2((z+l)/J-x) \,dx},\nonumber %
\end{eqnarray}
where $\hth_j$ are chosen so that $\hf_j$ satisfy (\ref{constr}).

The solution of (\ref{mi}) is not given
explicitly. The estimates are calculated by an iterative algorithm with
a starting set of function estimates $\hf_1^{[0]}$ and
$\hf_2^{[0]}$ that satisfy the constraints (\ref{constr}).
With the initial estimates, we compute $\hf_3^{[0]}$ from the third
equation at (\ref{mi}). Then we update
$\hf_j^{[k-1]}$ consecutively for $j=1,2,3$ and for $k \ge1$ by the
equations at (\ref{mi}) until convergence. Specifically,
we compute at the $k$th cycle ($k \ge1$) of the iteration
%
\begin{eqnarray}
\label{alg} %
\hf_1^{[k]}(x) &=&
\hth_1^{[k]} \cdot\frac{\hf_{w,1}(x)}{\int_{J_2(x)}
\hf_2^{[k-1]}(y) \hf_3^{[k-1]}(\mJ(x+y)) \,dy},\nonumber
\\
\hf_2^{[k]}(y) &=& \hth_2^{[k]} \cdot
\frac{\hf_{w,2}(y)}{\int_{J_1(y)}
\hf_1^{[k]}(x) \hf_3^{[k-1]}(\mJ(x+y)) \,dx},
\\
\hf_3^{[k]}(z) &=& \hth_3^{[k]} \cdot
\frac{\hf_{w,3}(z)}{\sum_{l=0}^{L(J)}\int_{J_{3l}(z)}\hf_1^{[k]}(x)
\hf_2^{[k]}((z+l)/J-x) \,dx},\nonumber %
\end{eqnarray}
where $\hth_j^{[k]}$ are chosen so that the resulting $\hf_j^{[k]}$
satisfy (\ref{constr}).

We note that the naive two-dimensional kernel density estimator is not
consistent near the boundary region,
which jeopardizes the properties of the solution
of the backfitting equation (\ref{mi}) at boundaries.
For a preliminary estimator $\hf$ of the joint density $f$, we take the
local linear estimation technique. The local linear estimator
$\hf$ we consider here is similar in spirit to the proposal of \citet{Che97}.
Let ${\mathbf a}(u,v;x,y)=(1, (u-x)/h_1, (v-y)/h_2)^\top$ and define
\[
\bA(x,y)=\int_S {\mathbf a}(u,v;x,y){\mathbf
a}(u,v;x,y)^\top h_1^{-1}h_2^{-1}K
\biggl(\frac{u-x}{h_1} \biggr) K \biggl(\frac{v-y}{h_2} \biggr) \,du \,dv,
\]
where $(h_1,h_2)$ is the bandwidth vector and $K$ is a symmetric
univariate probability density function. Also, define
\[
\hat\bb(x,y)=n^{-1}\sum_{i=1}^n
{\mathbf a}(X_i,Y_i;x,y) h_1^{-1}h_2^{-1}K
\biggl(\frac{X_i-x}{h_1} \biggr) K \biggl(\frac{Y_i-y}{h_2}
\biggr)W_i,
\]
where $W_i=1$ if $(X_i,Y_i)\in S$ and $0$ otherwise.
The local linear density estimator $\hf$ we consider in this paper is
defined by $\hat\eta_0$, where
$\hat\etam=(\hat\eta_0,\hat\eta_1,\hat\eta_2)$ is given by
%
\begin{equation}
\label{loclinest} \hat\etam(x,y) = \bA(x,y)^{-1}\hat\bb(x,y).
\end{equation}
It is alternatively defined as
\begin{eqnarray*}
\hat\etam(x,y)&=&\argmin_\etam\lim_{b_1, b_2 \rightarrow0}\int
_S \bigl[ \hf_{b_1,b_2}(u,v) - {\mathbf
a}(u,v;x,y)^\top\etam(x,y) \bigr]^2
\\
&&\hspace*{79pt}{} \times K \biggl(\frac{u-x}{h_1} \biggr) K \biggl(\frac{v-y}{h_2} \biggr)
\,du \,dv,
\end{eqnarray*}
where $\hf_{b_1,b_2}$ be the standard two-dimensional kernel density
estimator defined by
\[
\hf_{b_1,b_2}(x,y) = n^{-1}\sum_{i=1}^n
b_1^{-1}b_2^{-1} K \biggl(
\frac
{x-X_i}{b_1} \biggr) K \biggl(\frac{y-Y_i}{b_2} \biggr)W_i
\]
for a bandwidth vector $(b_1,b_2)$.

Before we close this section, we give two remarks.
One is that, instead of integrating the two-dimensional estimator $\hf$,
one may estimate $f_{w,j}$ directly from the data. In particular, one
may estimate $f_{w,j}$ by the one-dimensional kernel
density estimators
\begin{eqnarray*}
\tf_{w,1}(x) &=&n^{-1}h_1^{-1}\sum
_{i=1}^n K \biggl(\frac
{X_i-x}{h_1}
\biggr)W_i,
\\
\tf_{w,2}(y) &=&n^{-1}h_2^{-1}\sum
_{i=1}^n K \biggl(\frac
{Y_i-y}{h_2}
\biggr)W_i,
\\
\tf_{w,3}(z) &=& n^{-1}h_3^{-1}\sum
_{i=1}^n K \biggl(\frac
{m_J(X_i+Y_i)-z}{h_3}
\biggr)W_i.
\end{eqnarray*}
Our theory that we present in the next section is valid for this
alternative estimation procedure. The other thing
we would like to remark is that one may be also interested in an
extension of the model (\ref{model}) that
arises when one observes a covariate $\bU_i \in\R^d$ along with
$(X_i,Y_i)$. A natural extension of the model (\ref{model})
in this case is that the conditional density of $(X, Y )$ given $\bU=
{\mathbf u}$ has the form
$f(x,y|{\mathbf u})=f_1(x,{\mathbf u})f_2(y,{\mathbf u})f_3(m_J(x+y),{\mathbf u}),
(x,y) \in\cI$, where
the constraints (B1) now applies to $f_1(\cdot,\bz)$ and $f_2(\cdot
,\bz
)$ for each $\bz$. The method and theory for this
extended model are easy to derive from those we present here.

\section{Theoretical properties}\label{sec4}

Let $\cS$ denote the space of function tuples $\bg=(g_1,g_2,g_3)$
with square integrable univariate functions $g_j$ in the space $L_2[0,1]$.
Define nonlinear functionals $\cF_j$ for $1 \le j \le3$ on $\cS$ by
\begin{eqnarray*}
\cF_1(\bg) &=& 1-\int_{S_1} g_1(x) \,dx,
\\
\cF_2(\bg) &=& 1-\int_{S_2} g_2(y) \,dy,
\\
\cF_3(\bg) &=& \vartheta- \int_S
g_1(x)g_2(y)g_3\bigl(\mJ(x+y)\bigr) \,dx \,dy.
\end{eqnarray*}
Also, define nonlinear functionals $\cF_j$ for $4 \le j \le6$, now on
$\R^3 \times\cS$, by
\begin{eqnarray*}
\cF_4(\thetam,\bg) (x) &=& \int_{J_2(x)} \bigl[
\theta_1 f(x,y)-g_1(x)g_2(y)g_3
\bigl(\mJ(x+y)\bigr) \bigr] \,dy,
\\
\cF_5(\thetam,\bg) (y) &=& \int_{J_1(y)} \bigl[
\theta_2 f(x,y)-g_1(x)g_2(y)g_3
\bigl(\mJ(x+y)\bigr) \bigr] \,dx,
\\
\cF_6(\thetam,\bg) (z) &=& \sum_{l=0}^{L(J)}
\int_{J_{3l}(z)} \bigl[\theta_3 f\bigl(x,(z+l)/J-x
\bigr)\\
&&\hspace*{48pt}{}-g_1(x)g_2\bigl((z+l)/J-x\bigr) g_3(z)
\bigr] \,dx,
\end{eqnarray*}
where $\thetam=(\theta_1,\theta_2,\theta_3)^\top$.
Then we define a nonlinear operator $\cF\dvtx  \R^3 \times\cS\mapsto\R^3
\times\cS$ by
$\cF(\thetam,\bg)(x,y,z) = (\cF_1(\bg), \cF_2(\bg), \cF_3(\bg
), \cF
_4(\thetam,\bg)(x),
\cF_5(\thetam,\bg)(y), \break  \cF_6(\thetam,\bg)(z))^\top$.

Now, we define nonlinear functionals $\hcF_j$ for $1 \le j \le3$ on
$\cS$ and $\hcF_j$ for $4 \le j \le6$
on $\R^3 \times\cS$ as $\cF_j$ in the above, with
the joint density $f$ being replaced by its estimator $\hf$ and
$\vartheta$ by $\hat\vartheta$.
Let $\hcF\dvtx  \R^3 \times\cS\mapsto\R^3 \times\cS$ be the nonlinear
operator defined by
$\hcF(\thetam,\bg)(x,y,z) = (\hcF_1(\bg), \hcF_2(\bg), \hcF
_3(\bg), \hcF
_4(\thetam,\bg)(x),
\hcF_5(\thetam,\bg)(y), \break \hcF_6(\thetam,\bg)(z))^\top$.
Our estimators $\hat\bff=(\hf_1,\hf_2,\hf_3)$ along with $\hat
\thetam
=(\hth_1, \hth_2, \hth_3)$ are
given as the solution of the equation
%
\begin{equation}
\label{est-def} \hcF(\hat\thetam,\hat\bff)=\bzero.
\end{equation}
From the definition of the nonlinear operator $\cF$, we also get
$\cF(\bone,\bff)=\bzero$, where $\bone=(1,1,1)^\top$ and $\bff
=(f_1,f_2,f_3)^\top$ for the true component
functions $f_j$.

We consider a theoretical approximation of $\hat\bff$. Define a
nonlinear operator by
$\cG(\thetam,\bg) = \cF(\bone+\thetam, \bff\circ(\bone+\bg
))$, where
$\bg_1\circ\bg_2$ denotes the entry-wise multiplication of the two
function vectors $\bg_1$ and $\bg_2$.
Then $\cG(\bzero,\bzero)=\bzero$. Let $\cG'(\bd,\deltam)$ denote the
derivative of $\cG(\thetam,\bg)$ at
$(\thetam,\bg)=(\bzero,\bzero)$ to the direction $(\bd,\deltam)$.
We write $\bff_w(x,y,z)=(f_{w,1}(x),f_{w,2}(y),f_{w,3}(z))^\top$ and
$\hat\mum(x,y,z)=(\hat\mu_1(x), \hat\mu_2(y),\hat\mu_3(z))^\top
$, where
%
\begin{eqnarray}
\label{defmu} %
\hat\mu_1(x) &=& f_{w,1}(x)^{-1}
\int_{J_2(x)} \bigl[\hf (x,y)-f(x,y) \bigr] \,dy,\nonumber
\\
\hat\mu_2(y) &=& f_{w,2}(y)^{-1}\int
_{J_1(y)} \bigl[\hf (x,y)-f(x,y) \bigr] \,dx,
\nonumber
\\[-8pt]
\\[-8pt]
\nonumber
\hat\mu_3(z) &=& f_{w,3}(z)^{-1}\sum
_{l=0}^{L(J)} \int_{J_{3l}(z)} \bigl[\hf
\bigl(x,(z+l)/J-x\bigr)
\\
&&\hspace*{95pt}{} - f\bigl(x,(z+l)/J-x\bigr) \bigr] \,dx.\nonumber %
\end{eqnarray}
Let $\cG'^{-1}\dvtx \R^3 \times\cS\mapsto\R^3 \times\cS$ denote the
inverse of
$\cG'$, whose existence we will prove in the \hyperref[app]{Appendix}.
We define $\babff=(\bar f_1, \bar f_2, \bar f_3)$ along with
$\bathetam
=(\bar\theta_1, \bar\theta_2, \bar\theta_3)$ by
%
\begin{equation}
\label{appest} %
\pmatrix{\bathetam-
\bone
\cr
(\babff- \bff)/\bff} =
\cG'^{-1}\pmatrix{\bzero
\cr
-\bff_w \circ\hat\mum },
\end{equation}
where $\bg_1/\bg_2$ denotes the entrywise division of the function
$\bg
_1$ by $\bg_2$.

It can be seen that $\deltam= (\delta_1,\delta_2,\delta_3)^\top=
((\bar f_1-f_1)/f_1,(\bar f_2-f_2)/f_2, (\bar f_3-f_3)/f_3)^\top$
along with
$\bd=(d_1, d_2, d_3)^\top= (\bar\theta_1-1, \bar\theta_2-1,\bar
\theta
_3-1)^\top$
are given as the solution of the following system of integral equations:
%
\begin{eqnarray}
\label{backeq-app} %
\delta_1(x) &=& d_1 + \hat
\mu_1(x) - \int_{J_2(x)}\delta_2(y)
\frac
{f(x,y)}{f_{w,1}(x)} \,dy\nonumber\\
&&{} -\int_{J_2(x)}\delta_3\bigl(
\mJ(x+y)\bigr) \frac{f(x,y)}{f_{w,1}(x)} \,dy,\qquad x \in S_1,
\\
\delta_2(y) &=& d_2 + \hat\mu_2(y) - \int
_{J_1(y)}\delta_1(x) \frac
{f(x,y)}{f_{w,2}(y)} \,dx
\nonumber\\
&&{}-\int
_{J_1(y)}\delta_3\bigl(\mJ(x+y)\bigr)
\frac{f(x,y)}{f_{w,2}(y)} \,dx,\qquad y \in S_2,
\nonumber\\
\delta_3(z) &=& d_3 + \hat\mu_3(z) - \sum
_{l=0}^{L(J)}\int_{J_{3l}(z)}
\delta_1(x) \frac
{f(x,(z+l)/J-x)}{f_{w,3}(z)} \,dx
\nonumber\\
&&{} - \sum_{l=0}^{L(J)}\int
_{J_{3l}(z)}\delta_2\bigl((z+l)/J-x\bigr)
\frac{f(x,(z+l)/J-x)}{f_{w,3}(z)} \,dx,\qquad z \in S_3,\nonumber %
\end{eqnarray}
subject to the constraints
%
\begin{eqnarray}
\label{constr-app} %
0&=&\int_{S_1} f_1(x)
\delta_1(x) \,dx,\nonumber
\\
0&=&\int_{S_2} f_2(y)\delta_2(y) \,dy,
\\
0&=&\int_S f(x,y) \bigl[\delta_1(x)+
\delta_2(y)+\delta_3\bigl(\mJ (x+y)\bigr) \bigr] \,dx \,dy.\nonumber
\end{eqnarray}

In the following theorem, we show that the approximation of $\hat\bff$
by $\babff$ is good enough. In the theorem, we
assume that $\hf(x,y)-f(x,y)=O_p(\ve_n)$ uniformly on $S$ for some
nonnegative sequence $\{\ve_n\}$ that converges to zero as
$n$ tends to infinity. For the local linear estimator $\hf$ defined by
(\ref{loclinest}) with $h_1 \sim h_2 \sim n^{-1/5}$,
we have $\ve_n = n^{-3/10}\sqrt{\log n}$.
The theorem tells that the approximation errors of $\bar f_j$ for $\hf
_j$ are of order $O_p(n^{-3/5}\log n)$. In Theorem~\ref{th4} below,
we will show that $\bar f_j-f_j$ have magnitude of order
$O_p(n^{-2/5}\sqrt{\log n})$ uniformly on $S_j$. This means that
the first-order properties of $\hf_j$ are the same as those of $\bar f_j$.

\begin{theorem}\label{th3} Assume that the conditions of Theorem~\ref{theoident2} hold, and that
the joint density $f$ is bounded away from zero
and infinity on its support $S$ with continuous partial derivatives on
the interior of $S$. If
$\hf(x,y)-f(x,y)=O_p(\ve_n)$ uniformly for $(x,y) \in S$, then it
holds that
$|\hat\theta_j - \bar\theta_j| = O_p(\ve_n^2)$ and
$\sup_{u \in S_j}|\hf_j(u)-\bar f_j(u)| = O_p(\ve_n^2)$.
\end{theorem}

Next, we present the limit distribution of $(\babff- \bff)/\bff$. In
the next theorem, we assume that $h_1 \sim c_1 n^{-1/5}$
and $h_2 \sim c_2 n^{-1/5}$ for some constants $c_1, c_2>0$. For
such constants, define
%
\begin{equation}
\label{tfB} \tf^B(x,y) = \frac{1}{2}\int u^2
K(u) \,du \biggl[c_1^2\frac{\partial
^2}{\partial x^2}f(x,y)
+c_2^2\frac{\partial^2}{\partial y^2}f(x,y) \biggr].
\end{equation}
Also, define $\tmu_j^B$ for $j=1,2,3$ as $\hat\mu_j$ at (\ref{defmu})
with the local linear estimator $\hf$ being replaced by $\tf^B$. In the
\hyperref[app]{Appendix}, we will show that the asymptotic mean of
$(\bar f_j -f_j)/f_j$ equals $n^{-2/5}\beta_j$, where $\betam= (\beta
_1, \beta_2, \beta_3)$ is
the solution of the backfitting equation (\ref{backeq-app}) with $\hat
\mum$ being replaced by
$\tilde\mum^B$. Let $\tf^A$ denote the
centered version of the naive two-dimensional kernel density estimator.
Specifically,
%
\begin{eqnarray}
\label{tfA} &&\tf^A(x,y) = n^{-1}\sum
_{i=1}^n \bigl[K_{h_1}(X_i-x)K_{h_2}(Y_i-y)
\nonumber
\\[-8pt]
\\[-8pt]
\nonumber
&&\hspace*{66pt}\qquad{}- E \bigl( K_{h_1}(X_i-x)K_{h_2}(Y_i-y)
\bigr) \bigr].
\end{eqnarray}
Here and below, we write $K_h(u)=K(u/h)/h$. Define $\tmu_j^A$ for
$j=1,2,3$ as $\tmu_j^B$ with $\tf^A$ taking the role of
$\tf^B$. We will also show that the asymptotic variances of $(\bar
f_j-f_j)/f_j$ equal those of $\tmu_j^A$, respectively, and
that they are given by $n^{-4/5}\sigma_j^2$, where
\begin{eqnarray*}
\sigma_1^2(x) &=& c_1^{-1}
f_{w,1}(x)^{-1}\int K^2(u) \,du,
\\
\sigma_2^2(y) &=& c_2^{-1}
f_{w,2}(y)^{-1}\int K^2(u) \,du,
\\
\sigma_3^2(z) &=& c_2^{-1}
f_{w,3}(z)^{-1}\int\bigl[K*K(u)\bigr] \bigl[K*K(c_1
u/c_2)\bigr] \,du
\\
&=& c_1^{-1} f_{w,3}(z)^{-1}\int
\bigl[K*K(u)\bigr] \bigl[K*K(c_2 u/c_1)\bigr] \,du,
\end{eqnarray*}
where $K*K$ denotes the two-fold convolution of the kernel $K$.

In the discussion of assumption (A6) in Section~\ref{sec2}, we note that (A6)
allows a finite number of jumps
in $I_j(u)$ for $j=1,2$ and $I_{3l}(u)$ as $u$ changes. These jump
points are actually those where the marginal
densities $f_{w,j}$ are discontinuous. At these discontinuity points,
the expression of the asymptotic distributions of
the estimators is complicate. For this reason, we consider only those
points in the partitions
$(a_{k-1}^j,a_k^j), 1 \le k \le L_j$, for the asymptotic distribution
of $\hf_j$,
where $a_k^j$ are the points that appear in assumption (A6). We denote
by $S_{j,c}$
the resulting subset of $S_j$ after deleting all $a_k^j, 1 \le k \le
L_j-1$. Note that $f_{w,j}$ is continuous
on $S_{j,c}$ due to (A6). In the theorem below, we also denote by
$S_j^\mathrm{ o}$ the interiors of $S_j$,
$j=1, 2, 3$.

For the limit distribution of $\hf_j$, we put an additional condition
on the support set. To state the condition, let $J_2^\mathrm{ o}(u_1;h_2)$
be a subset of $J_2(u_1)$ such that
$v \in J_2^\mathrm{ o}(u_1;h_2)$ if and only if $v-h_2 t \in J_2(u_1)$ for
all $t \in[-1,1]$.
The set $J_2^\mathrm{ o}(u_1;h_2)$ is inside $J_2(u_1)$ at a depth $h_2$.
In the following assumption, $a_k^j$
and $\kappa$ are the points and the function that appear in assumption (A6).
\begin{longlist}[(A7)]
\item[(A7)] There exist constants $C>0$ and $\alpha>1/2$ such that the
following statements hold:
(i) for any sequence of positive numbers $\varepsilon_n$,
$J_2^\mathrm{ o}(u_1;C \varepsilon_n^\alpha) \subset J_2(u_2)$ for all $u_1,
u_2 \in(a_{k-1}^1,a_k^1)\cap S_1$
with $|u_1-u_2| \le\varepsilon_n$, $1 \le k \le L_1$;
$J_1^\mathrm{ o}(u_1;\break C \varepsilon_n^\alpha) \subset J_1(u_2)$ for all $u_1,
u_2 \in(a_{k-1}^2,a_k^2)\cap S_2$
with $|u_1-u_2| \le\varepsilon_n$, $1 \le k \le L_2$; (ii)
$\kappa(t) \le C|t|^\alpha$.
\end{longlist}

\begin{theorem}\label{th4}
Assume that \textup{(A7)} and the conditions of Theorem~\ref{th3} hold, and that the
joint density $f$ is twice partially continuously differentiable.
Let the kernel $K$ be supported on $[-1,1]$, symmetric
and Lipschitz continuous. Let the bandwidths $h_j$ satisfy $n^{1/5}h_j
\rightarrow c_j$ for some constants $c_j>0$.
Then, for fixed points $u_j \in S_j^\mathrm{ o} \cap S_{j,c}$, it holds that
$n^{2/5}(\bar f_j(u_j)-f_j(u_j))/f_j(u_j)$ are jointly asymptotically
normal with mean $(\beta_j(u_j)\dvtx 1 \le j \le3)$
and variance $\operatorname{ diag}(\sigma^2_j(u_j)\dvtx 1 \le j \le3)$. Furthermore,
$(\bar f_j(u_j)-f_j(u_j))/f_j(u_j)=\break O_p(n^{-2/5}
\sqrt{\log n})$
uniformly for $u_j \in S_j$.
\end{theorem}

\begin{remark}
In the case where the third component function $f_3$ is constant, that
is, there is no periodic component,
the above theorem continue to hold for the component $f_1$ and $f_2$
without those conditions that
pertain to the set $S_3$ and the function $f_3$.
\end{remark}

\section{Numerical properties}\label{sec5}

\subsection{Simulation studies}

We considered two densities on $\cI=\{(x,y)\dvtx 0 \le x,y \le1, x+y
\le1\}$. Model 1 has
the components $f_1\equiv f_2 \equiv1$ on $[0,1]$, and $f_3(u) = c_1
(\sin(2\pi u)+3/2), u \in[0,1]$,
where $c_1>0$ is chosen to make $f(x,y)=f_1(x)f_2(y)f_3(\mJ(x+y))$ be a
density on $\cI$.
Model 2 has $f_1(u)=3/2-u, f_2(u)=5/4-3u^2/4$ and
$f_3(u)=c_2(u^3-3u^2/2 +u/2 +1/2)$ for some constant $c_2>0$.
We took $J=2$. We computed our estimates on a grid of bandwidth choice
$h_1=h_2$. For model 1, we took
$\{0.070+0.001\times j\dvtx 0 \le j \le30\}$ in the range $[0.070,0.100]$,
and for model 2 we chose
$\{0.40+0.02\times j\dvtx 0 \le j \le20\}$ in the range $[0.40,0.80]$. In
both cases, the ranges covered
the optimal bandwidths.
We obtained $\mathrm{ MISE}_j = E \int_0^1 [\hat f_j(u) - f_j(u)]^2 \,du$,
$\mathrm{ ISB}_j =
\int_0^1 [E \hat f_j(u) - f_j(u)]^2 \,du$ and $\mathrm{ IV}_j = E\int_0^1[\hat f_j(u)- E \hat f_j(u)]^2 \,du$,
for $1 \le j \le3$, based on 100 pseudo samples.
The sample sizes were $n=400$ and $1000$, but only the results for
$n=400$ are reported since
the lessons are the same.

\begin{figure}

\includegraphics{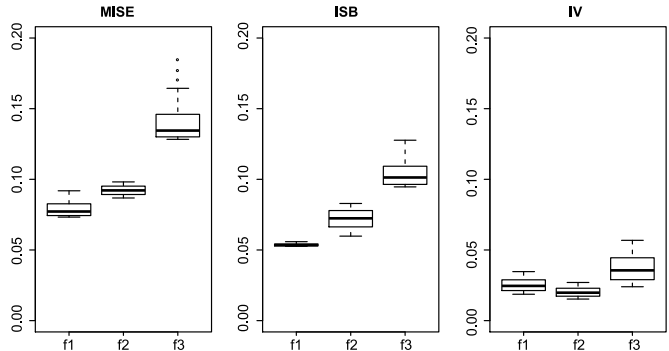}

\caption{Boxplots for the values of MISE, ISB and IV of our
estimates $f_j$ computed using various bandwidth choices (model 1, $n=400$).}
\label{fig:boxplot400}
\end{figure}

\begin{figure}[b]

\includegraphics{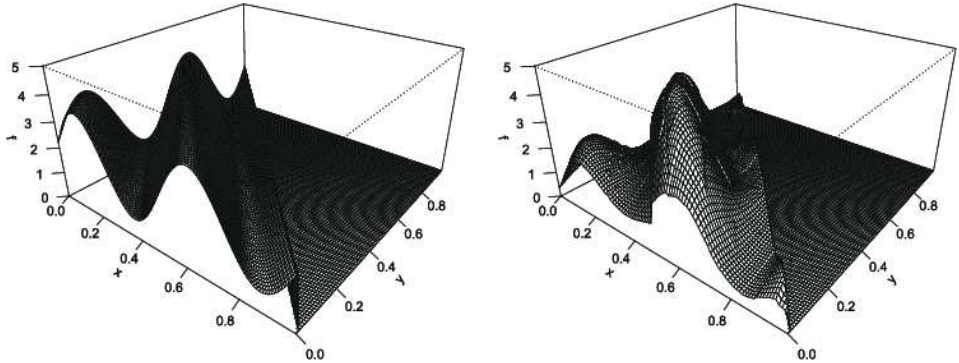}

\caption{The true density (left) and our estimated
two-dimensional density function (right)
computed from the pseudo sample that gives the median performance in
terms of ISE, for model 1 and $n=400$.}
\label{fig:2Dest}
\end{figure}

Figure~\ref{fig:boxplot400} is for model~1. It shows the boxplots of the
values of $\mathrm{MISE}_j$, $\mathrm{ISB}_j$ and $\mathrm{IV}_j$ computed using the bandwidths on
the grid specified above, and
thus gives some indication of how sensitive our estimators are to the
choice of bandwidth.
The bandwidth that gave the minimal value of $\mathrm{ MISE}_1+\mathrm{
MISE}_2+\mathrm{ MISE}_3$ was
$h_1=h_2=0.089$ in model~1, and $h_1=h_2=0.64$ in model~2, for the
sample size $n=400$.
The values of $\mathrm{MISE}_j$ along with $\mathrm{ISB}_j$ and $\mathrm{IV}_j$ for these optimal
bandwidths are reported in Table~\ref{tab1}.
Although our primary concern is the estimation of the component
functions, it is also of interest to see
how good the produced two-dimensional density estimator $\hat
f_1(x)\hat f_2(y)\hat f_3(\mJ(x+y))$ behaves. For this,
we include in the table the values of MISE, ISB and IV of the
two-dimensional estimates computed using the optimal bandwidth
$h_1=h_2=0.089$ in model~1, and $h_1=h_2=0.64$ in model~2.
For comparison, we also report the results for the two-dimensional
local linear
estimates defined at (\ref{loclinest}). For the local linear estimator,
we used its optimal choices $h_1=h_2=0.085$ in
model~1, and $h_1=h_2=0.48$ in model~2. We found that the initial local
linear estimates had a large portion of mass
outside $\cI$, and thus behaved very poorly if they were not re-scaled
to be integrated to one on $\cI$. The reported
values in Table~\ref{tab1} are for the adjusted local linear estimates. Overall,
our two-dimensional estimator has better performance
than the local linear estimator, especially in model~2. Figure~\ref{fig:2Dest} depicts the true density of model~1 and
our two-dimensional estimate that has the median performance in terms
of ISE.

\begin{table}
\caption{Mean integrated squared errors (MISE),
integrated squared biases (ISB) and integrated variance
(IV) of the estimators}\label{tab1}
\begin{tabular*}{\textwidth}{@{\extracolsep{\fill}}lcccccc@{}}
\hline
&&\multicolumn{3}{c}{\textbf{Component functions}}&\multicolumn{2}{c@{}}{\textbf{Joint density}}\\[-6pt]
&&\multicolumn{3}{c}{\hrulefill}&\multicolumn{2}{c@{}}{\hrulefill}\\
&&$\bolds{f_{1}}$&$\bolds{f_{2}}$&$\bolds{f_{3}}$&\textbf{Our est.} &\textbf{Local linear}\\
\hline
Model~1 &MISE &0.0756 &0.0937 &0.1283 &0.2493&0.2537 \\
&ISB &0.0528&0.0752 &0.0963 &0.1844&0.2199\\
&IV &0.0228 &0.0184 &0.0320&0.0649&0.0338 \\[3pt]
Model~2 &MISE &0.0124 &0.0057 &0.0130 &0.0475& 0.0624\\
&ISB &0.0120&0.0054 &0.0127 & 0.0469 &0.0607\\
&IV &0.0004 &0.0003 &0.0003& 0.0006&0.0017\\
\hline
\end{tabular*}
%
\end{table}

\subsection{Data examples}

The original data set we analyze in this section was collected between
the year 1990 to 2011
by the major global UK based nonlife insurance company RSA. The
dataset---and more details
about it---is publicly available via the Cass Business School web site
together with the paper
``Double Chain Ladder'' at the Cass knowledge site.
The observations were the incurred counts of large claims aggregated by months.
During the 264 months, 1516 large claims were made.
The dataset is provided in the form of a classical run-off triangle $\{
N_{kl}\dvtx 1 \le k, l \le264, k+l \le265\}$, where
$N_{kl}$ denotes the number of large claims incurred in the $k$th month
and reported in the $(k+l-1)$th month, that is, with
$(l-1)$ months delay.
Since the data are grouped monthly, we need pre-smoothing of the data
to apply the model (\ref{model})
that is based on data recorded over a continuous time scale. A natural
way of pre-smoothing is to perturb the
data by uniform random variables. Thus, we converted each claim $(k,l)$
on the two-dimensional discrete time scale
$\{(k,l)\dvtx 1 \le k, l \le264, k+l \le265\}$, into $(X,Y)$ on the
two-dimensional continuous time scale
$\cI=\{(x,y)\dvtx 0 \le x,y \le1, x+y \le1\}$, by
\[
X=\frac{k-1+U_1}{264},\qquad Y= \frac{l-1+U_2}{264},
\]
where $(U_1,U_2)$ is a two-dimensional uniform random variate on the
unit square $[0,1]^2$.
This gives a converted dataset $\{(X_i,Y_i)\dvtx 1 \le i \le1516\}$.
We applied to this dataset our method of estimating the structured
density $f$ of $(X,Y)$.

Since one month corresponds to an interval with length $1/264$ on the
$[0,1]$ scale, one year is equivalent to an interval
with length $12/264=1/22$ on the latter scale. We let the periodic
component $f_3(m_J(\cdot))$ in the model (\ref{model})
reflect a possible seasonal effect, so that we take one year in the
real time to be the period of the function.
This means that we let the periodic component $f_3(m_J(\cdot))$ have
$1/22$ as its period, and thus take $J=22$.
For the bandwidth, we took $h_1=h_2=0.01$.
The chosen bandwidth may be considered to be too small for the
estimation of $f_1$ and~$f_2$. However,
we took such a small bandwidth to detect possible seasonality. Note
that the bandwidth size $0.01$ corresponds to
$0.01\times12 \times22=2.64$ months. We found that even with this
small bandwidth
the estimated curve $\hf_3$ was nearly a constant function, which suggests
that the large claim data do not have a seasonal effect.

To see how well our method detects a possible seasonal effect in the
data, we augmented the dataset by adding a certain
level of seasonal effect as follows. We computed
\begin{eqnarray*}
N_{kl}' &=& 2 N_{kl} \qquad\mbox{if } k+l=12m
\mbox{ for some } m=1,2, \ldots,
\\
N_{kl}' &=& 3 N_{kl} \qquad\mbox{if } k+l=12m+1
\mbox{ for some } m=1,2, \ldots,
\\
N_{kl}' &=& 5 N_{kl}\qquad \mbox{if } k+l=12m+2
\mbox{ for some } m=0,1, \ldots,
\\
N_{kl}' &=& 3 N_{kl} \qquad\mbox{if } k+l=12m+3
\mbox{ for some } m=0,1, \ldots,
\\
N_{kl}' &=& N_{kl}\qquad \mbox{otherwise}.
\end{eqnarray*}
Since ($k+l-1$ modulo $12$) is the actual month of the claims reported,
the augmented dataset has added
claims in November, December, January and February.
The augmentation resulted in increasing the total number of claims
to 2606 from 1516. The increased counts of reported claims were 252
from 126 for November, 600 from 200 for December,
455 from 91 for January and 300 from 100 for February.

In our estimation procedure, the bandwidths $h_1$ and $h_2$ control the
smoothness of the local linear estimate $\hat f$ along the $x$- and
$y$-axis, respectively. Consequently, choosing small values for $h_1$
and $h_2$ would result in nonsmooth estimates of the functions $f_1$
and $f_2$, which we observed in the pilot study with $h_1=h_2=0.01$.
Nevertheless, in some cases setting these bandwidths to be small,
relative to the scales of $X$ and $Y$, might be preferred when one
needs to detect possible seasonality, as is the case with the current
dataset. In our dataset, the bandwidth size $1/264=0.0038$ on the scale
of $[0,1]$ corresponds to one month in real time. Thus, taking the
bandwidths to be $0.015$, for example, that corresponds to a period of
four months, forces the seasonal effect to almost vanish in the
estimate of $f_3$.

To achieve both aims of
producing smooth estimates of $f_1$ and $f_2$, and of detecting
possible seasonal effect, we applied to the augmented dataset
a two-stage procedure that is based on our estimation method described
in Section~\ref{methodology}. In the first stage, we got a local
linear estimate $\hat f$ with $h_1=h_2=0.01$, and found an estimate of
$f_3$ using the iteration scheme at (\ref{alg}).
In the second stage, we recomputed a local linear estimate $\hat f$
with larger bandwidths $h_1=h_2=0.05$, and found estimates
of $f_1$ and $f_2$ using only the first two updating equations at (\ref
{alg}) with $\hat f_3^{[k-1]}$ being replaced by
the estimate of $f_3$ obtained in the first stage.

\begin{figure}

\includegraphics{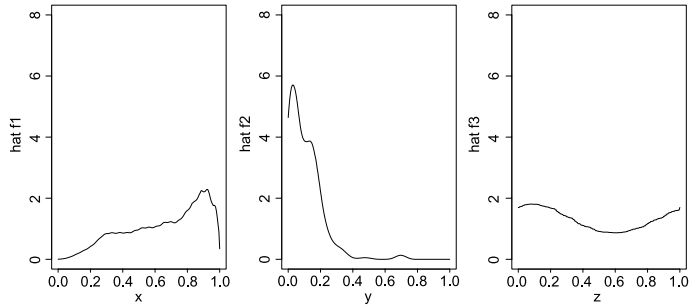}

\caption{Estimated curves $\hat f_j$ for the model (\protect\ref
{model}) obtained by applying the two-stage procedure
to the augmented large claim data.}
\label{fig:compest2}
\end{figure}

The results of applying this two-stage procedure to the augmented dataset
are presented in Figure~\ref{fig:compest2}.
Clearly, the seasonal effect of the augmented dataset was well
recovered in the estimate of $f_3$, and at the same time smooth estimates
of $f_1$ and $f_2$ were produced. The augmented data set indicate an
increased number of claims in the winter time. This is clearly
reflected in the estimated results, where the first part and the last
part of the estimated effect is higher than the rest of the curve.
Imagine the realistic situation that a nonlife insurer on the first
day of November has to produce budget expenses for the rest of the
year. The classical multiplicative methodology is not able to reflect
the two month perspective of such a budget. Therefore, considerable
work is being done manually in finance and actuarial departments of
nonlife insurance companies to correct for such effects.
With our new seasonal correction, costly manual procedures can be
replaced by cost saving automatic ones eventually benefitting the
prices all of us as end customers have to pay for insurance products.

\begin{figure}

\includegraphics{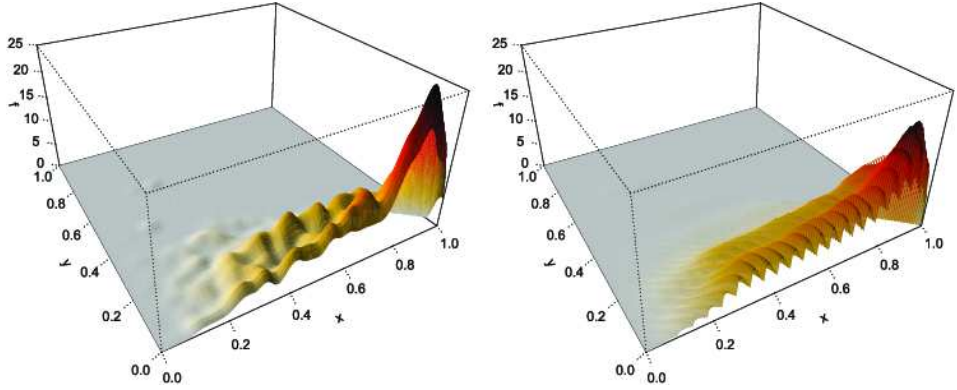}

\caption{Local linear joint density estimate (left) and our
estimate (right) for the model (\protect\ref{model})
obtained by applying the two-stage procedure to the augmented large
claim data.}
\label{fig:joint}
\end{figure}

Figure~\ref{fig:joint} depicts the resulting two-dimensional joint
density. Notice that this two-dimensional density is clearly
nonmultiplicative. The seasonal correction provides a visually
deviation from the multiplicative shape. Also, note that while this
two-dimensional density is nonmultiplicative, the nature of this
deviation is not immediately clear to the eye. Whether the deviation is
pure noise, a seasonal effect or some other effect is not easy to get
from the full two-dimensional graph of the
local linear density estimate which is also presented in Figure~\ref{fig:joint}. For the local linear estimate, we used $h_1=h_2=0.03$.
We tried other bandwidth choices such as $0.01$ and $0.05$, but found
that the smaller one gave too rough estimate and the
larger one produced too smooth a surface.
Our two-dimensional density estimate therefore illustrates why research
into structured densities on nontrivial supports is crucial to extract
information beyond the classical and simple multiplicative one.

\begin{appendix}\label{app}
\section*{Appendix}

\subsection{Proof of Theorem~\texorpdfstring{\protect\ref{theoident1}}{1}}

Suppose that $(g_{1},g_{2},g_{3}) $ is a tuple of functions that are
bounded away from zero and infinity with $\int_0^1 g_1(x) \,dx = \int_0^1 g_2(y) \,dy =1$ and
\[
f(x,y) = g_{1}(x) g_{2}(y) g_{3}\bigl(\mJ(x+y)
\bigr).
\]
Furthermore, we assume that $g_{1}$ and $g_{2}$ are differentiable on
$[0,1]$ and that $g_{3}$ is twice differentiable on $[0,1)$.
For $j \in\{1,2,3\}$ define $\mu_j = \log f_j - \log g_j$. By
assumption, we have
\[
\mu_1(x) + \mu_2(y) + \mu_3
\bigl(m_J(x+y)\bigr)=0.
\]

For $z \in[0,1)$, we choose $(x,y)$ in the interior of $\cI$ with
$m_J(x+y)=z$.
Then we have that
\[
0 = {\partial^2 \over\partial x\, \partial y} \bigl[\mu_1(x) + \mu_2(y)+
\mu _3\bigl(m_J(x+y)\bigr)\bigr] = \mu_3^{\prime\prime}
(z).
\]
Thus, $\mu_3$ is a linear function. Furthermore, we have that $\mu_3(0)
= \mu_3(1-)$.
This follows by noting that $\mu_3(0)= - \mu_1(x) - \mu_2(y) $ for
$(x,y) \in\cI$ with $m_J(x+y)=0$. Note that
$m_J(x+y)=0$ if and only if $x+y=l/J$ for some $l \ge1$, if $(x,y)$ is
in the interior of $\cI$.
After slightly
decreasing $x$ and $y$ to $x+\delta_x$ and $y+\delta_y$ with small
$\delta_x <0$, $\delta_y<0$, we have that
$\mu_3(1+J(\delta_x+\delta_y))= - \mu_1(x+\delta_x) - \mu
_2(y+\delta_y)$
since $m_J(x+y+\delta_x+\delta_y)=1+J(\delta_x+\delta_y)$.
Thus, $\mu_3(0) = \mu_3(1-)$ follows from continuity of $\mu_1$ and $
\mu_2$.
We conclude that $ \mu_3$ must be a constant function. Thus, $\mu_1(x)
+ \mu_2(y)$ is a constant function.

From assumption (A5), we get that $\mu_1(x)$ is constant on the
intervals $[x_j, x_{j+1}]$. Because the union of these intervals is
equal to $[0,1]$ we conclude that $\mu_1(x)$ is constant on $[0,1]$.
Using again (A5) we get that $\mu_2(y)$ is constant on $[0,1]$. Because
of the assumption that $\int_0^1 g_1(x) \,dx = \int_0^1 g_2(y) \,dy =1$
and $\int_0^1 f_1(x) \,dx = \int_0^1 f_2(y) \,dy =1$ we get that $f_1=
g_1$, $f_2= g_2$ and $f_3= g_3$. This completes the proof.

\subsection{Proof of Theorem~\texorpdfstring{\protect\ref{theoident2}}{2}}

We first argue that $\mu_1$, $\mu_2$ and $\mu_3$ are a.e. equal to
piecewise continuous functions on $(0,1)$,
with a finite number of pieces.
To see that $\mu_1$ is a.e. equal to a piecewise continuous function,
we note that
\[
\mu_1(x) = - \int_{I_2(x)} \bigl[
\mu_2(y) + \mu_3\bigl(m_J(x+y)\bigr) \bigr]
\,dy /\mes\bigl(I_2(x)\bigr)\qquad \mbox{a.e. } x\in(0,1).
\]
Here, because of (A3) and (A6), the right-hand side is a piecewise
continuous function. Thus, $\mu_1$ is a.e. equal to a piecewise
continuous function. In abuse of notation, we now denote the piecewise
continuous function by $\mu_1$.
By similar arguments, one sees that $\mu_2$, and $ \mu_3$ are piecewise
continuous functions (or more precisely a.e. equal to piecewise
continuous functions).
This implies that
%
\begin{equation}
\label{zerohelp1} \mu_1(x) + \mu_2(y) + \mu_3
\bigl(m_J(x+y)\bigr)=0
\end{equation}
for $(x,y,m_J(x+y)) \notin\{x_1,\ldots,x_{r_1}\} \times(0,1)^2 \cup
(0,1) \times\{y_1,\ldots,y_{r_2}\} \times(0,1) \cup(0,1)^2 \times\{
z_1,\ldots,z_{r_3}\}$ for some values $x_1,\ldots,x_{r_1}, y_1,\ldots,y_{r_2},
z_1,\ldots,z_{r_3} \in(0,1)$.

We now argue that $\mu_3$ is continuous on $[0,1)$. To see that $\mu_3$
is continuous at $z_0 \in[0,1)$,
we choose $(x_0,y_0)$ in the interior of $\cI$ such that $m_J(x_0+y_0)=z_0$.
This is possible because of assumption (A2). We can choose $x_0$ and
$y_0$ such that $\mu_1$ is continuous
at $x_0$ and $\mu_2$ is continuous at $y_0$.
Thus, we get from (\ref{zerohelp1}) that $\mu_3$ is continuous at
$z_0$. Similarly, one shows that
$\mu_1$ and $\mu_2$ are continuous functions on $[0,1]$. This gives that
%
\begin{equation}
\label{zerohelp2} \mu_1(x) + \mu_2(y) + \mu_3
\bigl(m_J(x+y)\bigr)=0
\end{equation}
for all $x,y \in(0,1)$.

For $z_0 \in[0,1)$, we choose $(x_0,y_0)$ in the interior of $\cI$
with $m_J(x_0+y_0)=z_0$.
Note that for $\delta_x$ and $\delta_y$ sufficiently small we get for
$z_0 \in(0,1)$ that
$m_J(x_0+ \delta_x + y_0+ \delta_y) = z_0 +J ( \delta_x + \delta_y)$.
This gives for $\delta_x$ and $\delta_y$ sufficiently small that
\[
\mu_1(x_0 + \delta_x ) +
\mu_2(y_0 + \delta_y ) + \mu_3
\bigl(z_0 +J ( \delta _x + \delta_y)
\bigr)=0.
\]
With $\delta_x$, $\delta_y'$ and $\delta_y$ sufficiently small, we
get that
\[
\mu_2(y_0 + \delta_y ) + \mu_3
\bigl(z_0 +J ( \delta_x + \delta_y)\bigr) =
\mu _2\bigl(y_0 + \delta_y'
\bigr) + \mu_3\bigl(z_0 +J \bigl( \delta_x +
\delta_y'\bigr)\bigr).
\]
With the special choice $\delta_x = -\delta_y$, this gives
\[
\mu_2(y_0 + \delta_y ) +
\mu_3(z_0 ) = \mu_2\bigl(y_0 +
\delta_y'\bigr) + \mu _3
\bigl(z_0 +J \bigl(\delta_y'-
\delta_y\bigr)\bigr).
\]
Let $\gamma$ be a function defined by $\gamma(u) = \mu_3(z_0 +J u) -
\mu
_3(z_0 )$.
From the last two equations taking $u=\delta_x+\delta_y$ and
$v=\delta
_y'-\delta_y$,
we get
\[
\gamma(u+v) = \gamma(u) + \gamma(v)
\]
for $u,v$ sufficiently small.
This implies that, with a constant $c_{z_0}$ depending on $z_0$ we have
$\gamma(u) = c_{z_0} u$ for $u$ sufficiently small; see Theorem~3 of
Guillot, Khare and
Rajaratnam (\citeyear{GuiKhaRaj13}).
Thus, we obtain $\mu_3(z) = a_{z_0} + b_{z_0} z$ with constants
$a_{z_0}$ and $b_{z_0}$ depending on $z_0$
for $z$ in a neighborhood $U_{z_0}$ of $z_0$. Because every interval
$[z',z'']$ with $0 < z' < z'' < 1$ can be covered
by the union of finitely many $U_{z}$'s we get that for each such
interval it holds that $\mu_3(z) = a_{z',z''}
+ b_{z',z''} z$ for $z \in[z',z'']$ with constants $a_{z',z''}$ and $
b_{z',z''}$ depending on the chosen interval
$[z',z'']$.

One can repeat the above arguments for $z_0=0$. Then we have that
$m_J(x_0+ \delta_x + y_0+ \delta_y) =
1 +J ( \delta_x + \delta_y)$ for $\delta_x + \delta_y < 0$ and
$m_J(x_0+ \delta_x + y_0+ \delta_y) =
J ( \delta_x + \delta_y)$ for $\delta_x + \delta_y > 0$. Arguing as
above with $\delta_x+\delta_y >0 $
and $\delta_y'- \delta_y> 0$ we get that $\mu_3(z) = a_{+} + b_{+} z$
for $z \in(0,z^{+}]$ for $z^{+}> 0 $
small enough with some constants $a_{+} $ and $ b_{+} $. Similarly, we
get by choosing $\delta_x+\delta_y <0 $
and $\delta_y'- \delta_y< 0$ that $\mu_3(z) = a_{-} + b_{-} z$ for $z
\in(z^{-},1)$ for $z^{-}<1 $
large enough with some constants $a_{-} $ and $ b_{-} $. Thus, we get that
$\mu_3(z) = a + b z$ for $z \in(0,1)$ with some constants $a $ and $ b$.

Furthermore, using continuity of $\mu_1$, $\mu_2$ and the relation
$\mu
_3(m_J(x+y)) = -\mu_1(x) - \mu_2(y) $ for $z=m_J(x+y)$ with $z$ in
$(1-\delta, 1)$ and $(0, \delta)$ with $\delta> 0$ small enough we get
that $\mu_3(0) = \mu_3(1-)$. Thus, we have $b=0$ and we conclude that
$\mu_3 $ is a constant function. This gives
\[
\mu_1(x) + \mu_2(y) = - a
\]
for all $(x,y) \in\cI$. Now arguing as in the proof of Theorem~\ref
{theoident1} we get that $f_1= g_1$, $f_2= g_2$ and $f_3= g_3$. This
completes the proof.

\subsection{Proof of Theorem~\texorpdfstring{\protect\ref{th3}}{3}}

Let $\cG'(\thetam,\bg)(\bd,\deltam)$ denote the derivative $\cG$,
defined in Section~\ref{sec4},
at $(\thetam,\bg)$ to the direction $(\bd,\deltam)$.
We note that we write $\cG'(\bzero,\bzero)(\bd,\deltam)$ simply as
$\cG'(\bd,\deltam)$ in Section~\ref{sec4}.
We use the sup-norm $\|(\bd,\deltam)\|_\infty$ as a metric in the space
$\R^3 \times\cS$, defined by
\[
\bigl\|(\bd,\deltam)\bigr\|_\infty= \max \Bigl\{|d_1|,
|d_2|, |d_3|, \sup_{u
\in
S_1}\bigl|
\delta_1(u)\bigr|, \sup_{u \in S_2}\bigl|\delta_2(u)\bigr|,
\sup_{u \in S_3}\bigl|\delta_3(u)\bigr| \Bigr\}.
\]
Define $\hcG(\thetam, \bg) = \hcF(\bone+\thetam, \bff\circ
(\bone+\bg
))$, where $\hcF$ is defined in Section~\ref{sec4},
and let $\hcG'(\thetam,\bg)$ denote
the derivative of $\hcG$ at $(\thetam,\bg)$. In the setting where
$\hf(x,y)-f(x,y) = O_p(\ve_n)$ uniformly for $(x,y) \in\cI$, we claim:
\begin{longlist}[(iii)]
\item[(i)] $\sup_{\|(\bd,\deltam)\|_\infty=1}\|\hcG'(\bzero
,\bzero)(\bd
,\deltam)
-\cG'(\bzero,\bzero)(\bd,\deltam)\|_\infty= O_p(\ve_n)$;
\item[(ii)] The operator $\cG'(\bzero,\bzero)$ is invertible and has
bounded inverse;
\item[(iii)] The operator $\hcG'$ is Lipschitz continuous with
probability tending to one, that is, there exists
constants $r, C>0$ such that, with probability tending to one,
\[
\sup_{\|(\bd,\deltam)\|_\infty=1}\bigl\|\hcG'(\thetam_1,
\bg_1) (\bd ,\deltam )-\hcG'(\thetam_2,
\bg_2) (\bd,\deltam)\bigr\|_\infty \le C\bigl \|(\thetam_1,
\bg_1)-(\thetam_2,\bg_2)\bigr\|_\infty
\]
for all $(\thetam_1,\bg_1), (\thetam_2,\bg_2) \in B_r(\bzero
,\bzero)$,
where $B_r(\thetam,\bg)$ is a ball with radius $r>0$ in
$\R^3 \times\cS$ centered at $(\thetam,\bg)$.
\end{longlist}

Theorem~\ref{th3} basically follows from the above (i)--(iii). To prove the
theorem using (i)--(iii), we note that claim (ii)
with the definitions of $\bathetam$ and $\babff$ at (\ref{appest}) gives
$\bathetam-\bone= O_p(\ve_n)$ and $(\babff-\bff)/\bff= O_p(\ve_n)$.
With (i) and (iii), this implies that
%
\begin{equation}
\label{lipsch} \sup_{\|(\bd,\deltam)\|_\infty=1}\bigl\|\hcG'\bigl(\bathetam-
\bone,(\babff -\bff )/\bff\bigr) (\bd,\deltam)-\cG'(\bzero,\bzero) (
\bd,\deltam) \bigr\| = O_p(\ve_n).
\end{equation}
Now, from (ii) it follows that there exists a constant $C>0$ such that
the map $\hcG'(\bathetam-\bone,(\babff-\bff)/\bff)$ is invertible
and $\|\hcG'(\bathetam-\bone,(\babff-\bff)/\bff)^{-1}(\bd
,\deltam)\|
_\infty\le C \|(\bd,\deltam)\|_\infty$
with probability tending to one. Also, (iii) is valid for all
$(\thetam_1,\bg_1),\break  (\thetam_2,\bg_2) \in B_{2r}(\bathetam-\bone
,(\babff-\bff)/\bff)$.
Then we can argue that the solution of the equation
$\hcG(\thetam,\bg)=\bzero$, which is $(\hthetam-\bone, (\hbff
-\bff)/\bff
)$, is within $C \alpha_n$ distance from
$(\bathetam-\bone,(\babff-\bff)/\bff)$, with probability tending to
one, where $C>0$ is a constant and
$\alpha_n = \|\hcG(\bathetam-\bone,(\babff-\bff)/\bff)\|_\infty
$. This follows
from an application of the Newton--Kantorovich theorem; see Deimling
(\citeyear{De85}) or \citet{YuParMam08} for
a statement of the theorem and related applications.
To compute $\alpha_n$, we note that
%
\begin{eqnarray}
\label{apperr} %
\hcG\bigl(\bathetam-\bone,(\babff-\bff)/\bff\bigr)&=&
\hcG(\bzero,\bzero) + \hcG '(\bzero,\bzero) \bigl(\bathetam-\bone,(
\babff-\bff)/\bff\bigr) + O_p\bigl(\ve_n^2
\bigr)
\nonumber
\\[-8pt]
\\[-8pt]
\nonumber
&=&\hcG(\bzero,\bzero) + \cG'(\bzero,\bzero) \bigl(\bathetam-\bone
,(\babff-\bff )/\bff\bigr) + O_p\bigl(\ve_n^2
\bigr). %
\end{eqnarray}
For the first equation of (\ref{apperr}), we have used (iii) and the
facts that $\bathetam-\bone= O_p(\ve_n)$ and $(\babff-\bff)/\bff=
O_p(\ve_n)$. The second equation of (\ref{apperr}) follows from the inequality
\[
\bigl\|\hcG'(\bzero,\bzero) (\bd,\deltam)-\cG'(\bzero,
\bzero) (\bd ,\deltam)\bigr\| _\infty\le C \sup_{x,y \in S}\bigl|
\hf(x,y) -f(x,y)\bigr| \cdot\bigl\|(\bd, \deltam)\bigr\|_\infty
\]
for some constant $C>0$. Now, $\hcG(\bzero,\bzero)=\hcF(\bone,\bff
)=(\bzero^\top,(\bff_w \circ\hat\mum)^\top)^\top$.
From the definition (\ref{appest}), we also get $\cG'(\bzero,\bzero
)(\bathetam-\bone,(\babff-\bff)/\bff)
= (\bzero^\top,-(\bff_w \circ\hat\mum)^\top)^\top$. This proves
$\alpha
_n=O_p(\ve_n^2)$, so that
$\|(\hthetam-\bathetam, (\hbff-\babff)/\bff)\|_\infty= O_p(\ve_n^2)$.

Claim (i) follows from the uniform convergence of $\hf$ to $f$ that is
assumed in the theorem:
$\sup_{(x,y)\in S}|\hf(x,y)-f(x,y)| = O_p(\ve_n)$. Below, we give the
proofs of claims (ii) and (iii).

\begin{pf*}{Proof of claim \normalfont{(ii)}}
For this claim, we first prove that the map $\cG'(\bzero,\bzero)$ is
one-to-one.
Suppose that $\cG'(\bzero,\bzero)(\bd,\deltam)=
\bzero$ for some $\bd=(d_1,d_2,d_3)^\top$ and $\deltam=(\delta
_1,\delta
_2,\delta_3)^\top$.
Then, by integrating the fourth component of $\cG'(\bzero,\bzero
)(\bd
,\deltam)$, we find that
\[
0=\int_S f(x,y)\bigl[\delta_1(x)+
\delta_2(y)+\delta_3\bigl(m_J(x+y)\bigr)
\bigr] \,dx \,dy = d_1 \int_S f(x,y) \,dx \,dy,
\]
where the first equation holds since the right-hand side equals, up to
sign change, the third component of
$\cG'(\bzero,\bzero)(\bd,\deltam)$.
Similarly, we get $d_2=d_3=0$. Now, from $\cG'(\bzero,\bzero)(\bzero
,\deltam)=\bzero$ we have
\begin{eqnarray*}
0&=&\int_{S_1 \times S_2 \times S_3}\bigl(\bzero^\top, \deltam
(x,y,z)^\top\bigr)\cG '(\bzero,\deltam) (x,y,z) \,dx \,dy \,dz
\\
&=&-\int_S f(x,y)\bigl[\delta_1(x)+
\delta_2(y)+\delta_3\bigl(m_J(x+y)\bigr)
\bigr]^2 \,dx \,dy.
\end{eqnarray*}
This implies
%
\begin{equation}
\label{ideneq} \delta_1(x)+\delta_2(y)+
\delta_3\bigl(m_J(x+y)\bigr) = 0 \qquad\mbox{a.e. on } S.
\end{equation}
Arguing as in the proof of Theorem~\ref{theoident2} using the last three equations
of
$\cG'(\bzero,\bzero)(\bzero,\deltam)=\bzero$, we obtain $\delta_j
\equiv0$ on $S_j$, $1 \le j \le3$.

Next, we prove that the map $\cG'(\bzero,\bzero)$ is onto. For a tuple
$(\bc,\etam)$ with $\bc=(c_1,c_2,c_3)^\top$ and
$\etam(x,y,z)=(\eta_1(x),\eta_2(y),\eta(z))^\top$, suppose that
$\langle(\bc,\etam), \break  \cG'(\bzero,\bzero)(\bd,\deltam)\rangle=0$
for all $(\bd, \deltam) \in\R^3 \times\cS$. This implies
%
\begin{eqnarray}
\label{onto} %
0&=&\int_S f(x,y)
\eta_1(x) \,dx \,dy,\nonumber
\\
0&=&\int_S f(x,y)\eta_2(y) \,dx \,dy,\nonumber
\\
0&=&\int_S f(x,y)\eta_3\bigl(m_J(x+y)
\bigr) \,dx \,dy,
\nonumber
\\
0&=&\int_{J_2(x)}f(x,y)\bigl[\eta_1(x)+
\eta_2(y)+\eta_3\bigl(m_J(x+y)\bigr)\bigr] \,dy
\nonumber\\
&&{}+ c_1f_1(x) + c_3f_{w,1}(x),
\\
0&=&\int_{J_1(y)}f(x,y)\bigl[\eta_1(x)+
\eta_2(y)+\eta_3\bigl(m_J(x+y)\bigr)\bigr] \,dx
\nonumber\\
&&{}+ c_2f_2(y) + c_3f_{w,2}(y),\nonumber
\\
\qquad0&=&\sum_{l=0}^{L(J)}\int_{J_{3l}(z)}
f\bigl(x,(z+l)/J-x\bigr)\bigl[\eta_1(x)+\eta _2
\bigl((z+l)/J-x\bigr)+ \eta_3(z)\bigr] \,dx\nonumber\\
&&{}+ c_3
f_{w,3}(z).\nonumber %
\end{eqnarray}
From the first three equations of (\ref{onto}), we get
$c_1+\vartheta c_3=0$ by integrating the fourth equation.
Similarly, we obtain $c_2+\vartheta c_3 =0$ and $c_3=0$ by integrating
the fifth and
the sixth equations. This establishes $c_1=c_2=c_3=0$. Putting back
these constant values to (\ref{onto}),
multiplying $\eta_1(x), \eta_2(y)$ and $\eta_3(z)$ to the right-hand
sides of the fourth, fifth and sixth equations, respectively,
and then integrating them give
\[
\int_S f(x,y)\bigl[\eta_1(x)+
\eta_2(y)+\eta_3\bigl(m_J(x+y)\bigr)
\bigr]^2 \,dx \,dy =0.
\]
Going through the arguments in the proof of $\cG'(\bzero,\bzero)$ being
one-to-one and now using the first two equations of
(\ref{onto}) give $\eta_1 = \eta_2 = \eta_3 \equiv0$. Note that the
first two equations can be written as
$\int_{S_1} f_{w,1}(x)\eta_1(x) \,dx = 0$ and $\int_{S_2}
f_{w,2}(y)\eta
_2(y) \,dy = 0$, and thus in the latter proof
$f_{w,j}$ for $j=1,2$ take the roles of $f_j$ in the former proof. The
foregoing arguments show that
$(\bzero,\bzero)$ is the only tuple that is perpendicular to the range
space of $\cG'(\bzero,\bzero)$,
which implies that $\cG'(\bzero,\bzero)$ is onto.

To verify that the inverse map $\cG'(\bzero,\bzero)^{-1}$ is bounded,
it suffices to prove that
the bijective linear operator $\cG'(\bzero,\bzero)$
is bounded, owing to the bounded inverse theorem. Indeed, it holds that
there exists a constant $C>0$ such that
$\|\cG'(\bzero,\bzero)(\bd,\deltam)\|_\infty\le C \|(\bd,\deltam
)\|
_\infty$. This completes the proof of claim~(ii).
\end{pf*}

\begin{pf*}{Proof of claim \normalfont{(iii)}} We first note that
$\hcG'(\thetam_1,\bg_1)(\bd,\deltam)- \hcG'(\thetam_2,\break\bg
_2) (\bd,\deltam)
=\cG'(\thetam_1,\bg_1)(\bd,\deltam)- \cG'(\thetam_2, \bg_2)(\bd
,\deltam)$.
From this, we get that, for each given $r>0$,
\[
\bigl\|\hcG'(\thetam_1,\bg_1) (\bd,\deltam) -
\hcG'(\thetam_2,\bg _2) (\bd ,\deltam)
\bigr\|_\infty\le 6 (1+r) \max_{1 \le j \le3}\sup
_{u \in S_j}f_{w,j}(u)\|\bg_2-\bg _1
\| _\infty
\]
for all $(\thetam_1,\bg_1), (\thetam_2,\bg_2) \in B_r(\bzero
,\bzero)$
and for all $(\bd,\deltam)$ with $\|(\bd,\deltam)\|_\infty=1$. For
this, we used
the inequality
\begin{eqnarray*}
&&\sup_{(x,y,z) \in S_1 \times S_2 \times S_3}\bigl|\kappa(x,y,z;\bg _2,\deltam )-
\kappa(x,y,z;\bg_1,\deltam)\bigr| \\
&&\qquad\le3 \|\deltam\|_\infty \bigl(2+\|
\bg_1\|_\infty+ \|\bg_2\|_\infty\bigr)\|
\bg_2-\bg_1\|_\infty.
\end{eqnarray*}
This completes the proof of (iii).
\end{pf*}

\subsection{Proof of Theorem~\texorpdfstring{\protect\ref{th4}}{4}}

Let $\hf^A(x,y)$ be the first entry of $\hat\etam^A(x,y)$, where
$\hat
\etam^A$ is defined
as $\hat\etam$ at (\ref{loclinest}) with $\hat\bb$ being replaced by
$\hat\bb- E\hat\bb$. Likewise, define $\hf^B(x,y)$
with $\hat\bb(x,y)$ being replaced by $E \hat\bb(x,y) - (f(x,y), h_1\,
\partial f(x,y)/\break \partial x,  h_2\, \partial f(x,y)/\partial y)^\top$.
Then $\hf(x,y)=f(x,y)+\hf^A(x,y)+\hf^B(x,y)$. Define
$\hat\mum^A$ and $\hat\mum^B$ as $\hat\mum$ at (\ref{defmu})
with $\hf
-f$ being replaced by $\hf^A$ and $\hf^B$,
respectively, and $\babff^s/\bff=(\bar f_1^s/f_1, \bar f_2^s/f_2,
\bar
f_3^s/f_3)$ along with
$\bar\thetam^s-\bone=(\bar\theta_1^s-1, \bar\theta_2^s-1,\bar
\theta
_3^s-1)$ for $s=A$ and $B$
as the solution of the backfitting equation (\ref{backeq-app}) with
$\hat\mum$ being replaced by $\hat\mum^s$,
subject to the constraints (\ref{constr-app}). Since the backfitting
equation (\ref{backeq-app}) is linear in $\hat\mum$,
we get that $\babff= \bff+ \babff^A + \babff^B$ and $\bar\thetam
=\bar
\thetam^A-\bone+ \bar\thetam^B$.

For simplicity, write the backfitting equation (\ref{backeq-app}) as
$\deltam=\bd+ \hat\mum- \bT\deltam$ with an appropriate
definition of the linear operator $\bT$. From the definitions of
$\babff
^A$ and $\bar\thetam^A$, we have
$\babff^A/\bff= \bar\thetam^A-\bone+ \hat\mum^A - \bT(\babff
^A/\bff
)$. From Lemma~\ref{lem1} below, we obtain
\[
\babff^A/\bff-\hat\mum^A= \bar\thetam^A-
\bone- \bT\bigl(\babff ^A/\bff-\hat \mum^A\bigr) +
o_p\bigl(n^{-2/5}\bigr)
\]
uniformly on $S_1 \times S_2 \times S_3$.
This implies $\babff^A/\bff-\hat\mum^A =o_p(n^{-2/5})$ uniformly on
$S_1 \times S_2 \times S_3$ and
$\bar\thetam^A-\bone= o_p(n^{-2/5})$.

Now, for the deterministic part $\babff^B$, recall the definitions of
$\tf^B$ and $\tilde\mum^B$ at (\ref{tfB})
and thereafter, respectively.
Let ${\mathbf r}_n = \hat\mum^B - n^{-2/5}\tilde\mum^B$. According to
Lemma~\ref{lem1}, ${\mathbf r}_n=o(n^{-2/5})$
on $S_1'\times S_2'\times S_3'$, where $S_j'$ is a subset of $S_j$ with
the property that $\mes(S_j-S_j')=O(n^{-1/5})$.
We also get ${\mathbf r}_n=O(n^{-2/5})$ on $S_1 \times S_2 \times S_3$.
This implies $\bT({\mathbf r}_n) = o(n^{-2/5})$, so that
\[
\babff^B/\bff-{\mathbf r}_n= \bar\thetam^B-
\bone+ n^{-2/5}\tilde\mum ^B - \bT\bigl(\babff^B/
\bff-{\mathbf r}_n\bigr) + o_p\bigl(n^{-2/5}
\bigr)
\]
uniformly on $S_1 \times S_2 \times S_3$.
Thus, $(\babff^B/\bff, \bar\thetam^B-\bone)$ equals the solution
of the
backfitting equation $\deltam= \bd+
n^{-2/5}\tilde\mum^B -\bT\deltam$, up to an additive term whose $j$th
component has a magnitude of an order
$o(n^{-2/5})$ on $S_j'$ and $O(n^{-2/5})$ on the whole set $S_j$.

The asymptotic distribution of $ ((\bar
f_j(u_j)-f_j(u_j))/f_j(u_j)\dvtx 1 \le j \le3 )$
for fixed $u_j \in S_{j,c}\cap S_j^\mathrm{ o}$
is then readily obtained from the above results. The asymptotic mean is
given as the solution
$(\delta_j(u_j)\dvtx 1 \le j \le3)$ of the backfitting equation (\ref
{backeq-app})
with $\hat\mu_j$ being replaced by $n^{-2/5}\tmu_j^B$, subject to the
constraint (\ref{constr-app}). The asymptotic variances
are derived from those of $\tmu_j^A$, where
\begin{eqnarray*}
\tmu_1^A(x)&=&f_{w,1}(x)^{-1}\int
_{J_2(x)}\tf^A(x,y) \,dy,
\\
\tmu_2^A(y)&=&f_{w,2}(y)^{-1}\int
_{J_1(y)}\tf^A(x,y) \,dx,
\\
\tmu_3^A(z)&=&f_{w,3}(z)^{-1}\sum
_{l=0}^{L(J)}\int_{J_{3l}(z)}\tf
^A\bigl(x,(z+l)/J-x\bigr) \,dx
\end{eqnarray*}
and $\tf^A(x,y) = n^{-1}\sum_{i=1}^n
[K_{h_1}(X_i-x)K_{h_2}(Y_i-y)W_i -
E ( K_{h_1}(X_i-x)K_{h_2}(Y_i-y)W_i ) ]$.
This is due to (\ref{tmu1-1}), (\ref{tmu2}) and the corresponding property
for $\hat\mu_3^A$ in the proof of Lemma~\ref{lem2} below.

To compute $\operatorname{ var}(\tmu_1^A(u_1))$, we note that, due to the
assumption (A7) and thus from Lemma~\ref{lem1},
we may find constants $C>0$ and $\alpha>1/2$ such that $J_2^\mathrm{
o}(u;Ch_1^\alpha+h_2) \subset
J_2^\mathrm{ o}(u_1;h_2)$ for all\vadjust{\goodbreak} $u$ with $|u-u_1|\le h_1$, if $n$ is
sufficiently large.
Note that $J_2^\mathrm{ o}(u;Ch_1^\alpha+h_2)$ is inside $J_2^\mathrm{ o}(u;h_2)$
at a depth $Ch_1^\alpha$.
Then it can be shown that, for all $(u,v)$ with $|u-u_1|\le h_1$ and $v
\in J_2^\mathrm{ o}(u;Ch_1^\alpha+h_2)$,
the set $\{(v-y)/h_2\dvtx y \in J_2(u_1)\}$ covers the interval $[-1,1]$,
the support of the kernel $K$.
This implies that
$K_{h_1}(u-u_1)\nu(u_1,v) = K_{h_1}(u-u_1)$ for all $(u,v)$ with
$|u-u_1|\le h_1$ and $v \in J_2^\mathrm{ o}(u;Ch_1^\alpha+h_2)$,
where
$\nu(u_1,v)=\int_{J_2(u_1)}K_{h_2}(v-y) \,dy$.
Using this and the fact that the Lebesgue measure of the set difference
$J_2(u)-J_2^\mathrm{ o}(u;Ch_1^\alpha+h_2)$
has a magnitude of order $n^{-\min\{1,\alpha\}/5}$, we get
\begin{eqnarray*}
&&\operatorname{ var}\bigl(\tmu_1^A(u_1)
\bigr)\\[-2pt]
&&\qquad=f_{w,1}(u_1)^{-2}n^{-1}h_1^{-1}
\int_S \frac{1}{h_1} K \biggl(\frac{u-u_1}{h_1}
\biggr)^2 \nu(u_1,v)^2 f(u,v) \,du \,dv + O
\bigl(n^{-1}\bigr)
\\[-2pt]
&&\qquad=f_{w,1}(u_1)^{-2}n^{-1}h_1^{-1}
\int_{|u-u_1|\le h_1}\int_{J_2^\mathrm{
o}(u;Ch_1^\alpha+h_2)} \frac{1}{h_1} K
\biggl(\frac{u-u_1}{h_1} \biggr)^2 \nu(u_1,v)^2
\\[-2pt]
&&\hspace*{192pt}\qquad\quad{} \times f(u,v) \,dv \,du \\[-2pt]
&&\qquad\quad{}+ o\bigl(n^{-1}h^{-1}\bigr)
\\[-2pt]
&&\qquad=f_{w,1}(u_1)^{-2}n^{-1}h_1^{-1}
\int_S \frac{1}{h_1} K \biggl(\frac{u-u_1}{h_1}
\biggr)^2 f(u,v) \,du \,dv + o\bigl(n^{-1}h^{-1}\bigr)
\\[-2pt]
&&\qquad= n^{-1}h_1^{-1} f_{w,1}(u_1)^{-1}
\int K^2(u) \,du + o\bigl(n^{-1}h^{-1}\bigr).
\end{eqnarray*}
The last equation holds since $u_1 \in S_{1,c}$, so that $f_{w,1}$ is
continuous at $u_1$,
and it is a fixed point in the interior of $S_1$.
Similarly, we obtain
\[
\operatorname{ var}\bigl(\tmu_2^A(u_2)\bigr) =
n^{-1}h_2^{-1} f_{w,2}(u_2)^{-1}
\int K^2(u) \,du + o\bigl(n^{-1}h^{-1}\bigr).
\]

The calculation of the asymptotic variance of $\tmu_3^A(u_3)$ is more
involved than those of $\operatorname{ var}(\tmu_j^A(u_j))$
for $j=1,2$. For this, we observe that, if $l \neq l'$, then for any
given $z \in[0,1]$ and $(u,v) \in\cI$ we have
\begin{eqnarray*}
&&\pi_{l,l'}\bigl(z,u,v,x,x'\bigr)
\\
&&\qquad \equiv K_{h_1}(u-x)K_{h_2} \biggl(v-\frac{z+l}{J}+x
\biggr)K_{h_1}\bigl(u-x'\bigr)K_{h_2} \biggl(v-
\frac{z+l'}{J}+x' \biggr)\\
&&\qquad =0
\end{eqnarray*}
for all $x, x'$ except the case $(z+l)/J-x = (z+l')/J-x'$, if $n$ is
sufficiently large. This implies that
\begin{eqnarray*}
&&\hspace*{-4pt}\operatorname{ var}\bigl(\tmu_3^A(u_3)
\bigr)\\
&&\hspace*{-8pt}\qquad=f_{w,3}(u_3)^{-2}n^{-1}\sum
_{l=0}^{L(J)} \int_{J_{3l}(u_3)}
\int_{J_{3l}(u_3)} \int_S \pi_{l}
\bigl(u_3,u,v,x,x'\bigr)f(u,v) \,du \,dv \,dx
\,dx'
\\
&&\hspace*{-8pt}\qquad\quad{} + O\bigl(n^{-1}\bigr),
\end{eqnarray*}
where $\pi_l = \pi_{l,l}$.
From Lemma~\ref{lem1} again, we may find constants $C>0$ and $\alpha>1/2$ such
that $J_2^\mathrm{ o}(x;Ch_1^\alpha+h_2) \subset
J_2^\mathrm{ o}(u;h_2)$ for all $x, u \in(a_{k-1}^1,a_k^1)\cap S_1$ with
$|u-x|\le h_1$,
$1 \le k \le L_1$.
Define a subset $J_{3l}'(u_3)$ of $[0,1]$ such that
$x \in J_{3l}'(u_3)$ if and only if $x \in J_{3l}(u_3+J(h_2+Ch_1^\alpha
)t)$ for all $t \in[-1,1]$.
Then, for a given $u \in S_{1,c}$, it follows that
\[
[-1,1] \subset \biggl\{ \frac{v-(u_3+l)/J +x}{h_2}\dvtx v \in J_2(u) \biggr
\}
\]
for all $x \in J_{3l}'(u_3)$ such that $|x-u|\le h_1$ and $x$ lies in
the same partition $(a_{k-1}^1,a_k^1)$ as $u$.
This holds since $x \in J_{3l}(z)$ implies
$(z+l)/J-x \in J_2(x)$. This entails that, for $x \in J_{3l}'(u_3) \cap
S_{1,c}^\mathrm{ o}(h_1)$,
\begin{eqnarray*}
&&\int_S \pi_l\bigl(u_3,u,v,x,x'
\bigr) \,du \,dv
\\
&&\qquad = \int_{[-1,1]^2} K(t)K(s)h_1^{-1}K
\biggl(t+\frac
{x-x'}{h_1} \biggr) h_2^{-1}K \biggl(s+
\frac{x'-x}{h_2} \biggr)\,dt \,ds
\\
&&\qquad =(K*K)_{h_1}\bigl(x-x'\bigr) (K*K)_{h_2}
\bigl(x-x'\bigr),
\end{eqnarray*}
where $K*K$ denotes the convolution of $K$ defined by $K*K(u)=\int
K(t)K(t+u) \,dt$.
Here and below, $S_{j,c}^\mathrm{ o}(h)$ for a small number $h>0$ denotes
the set of $x \in S_{j,c}$ such that
$x+ht$ belongs to $S_{j,c}$ for all $t \in[-1,1]$.

Because of the assumption (A7)
and the fact that $u_3$ is a fixed point in $S_{3,c}$,
we get that $\sum_{l=0}^{L(J)}\mes[J_{3l}(u_3) \triangle
J_{3l}'(u_3)]$ is of order $o(1)$.
This and the foregoing arguments give
\begin{eqnarray*}
&&\operatorname{ var}\bigl(\tmu_3^A(u_3)
\bigr)\\
&&\qquad=f_{w,3}(u_3)^{-2}n^{-1}\sum
_{l=0}^{L(J)} \int_{J_{3l}(u_3)}
\int_{J_{3l}'(u_3)\cap S_{1,c}^\mathrm{ o}(h_1)} \int_S \pi _{l}
\bigl(u_3,u,v,x,x'\bigr) \,du \,dv
\\
&&\hspace*{196pt}\qquad\quad{} \times f \biggl(x,\frac{u_3+l}{J}-x \biggr)
 \,dx \,dx'\\
 &&\qquad\quad{}+ o
\bigl(n^{-4/5}\bigr)
\\
&&\qquad=f_{w,3}(u_3)^{-2}n^{-1}\sum
_{l=0}^{L(J)} \int_{J_{3l}(u_3)} \int
_{J_{3l}(u_3)} (K*K)_{h_1}\bigl(x-x'\bigr)
(K*K)_{h_2}\bigl(x-x'\bigr)
\\
&&\hspace*{148pt}\qquad\quad{} \times f \biggl(x,\frac{u_3+l}{J}-x \biggr)
 \,dx \,dx'\\
 &&\qquad\quad{}+ o
\bigl(n^{-4/5}\bigr).
\end{eqnarray*}
Let $J_{3l}^\mathrm{ o}(u_3;2h_1)$ denote a subset of $J_{3l}(u_3)$ such
that $x \in J_{3l}^\mathrm{ o}(u_3;2h_1)$
if and only if $x-2h_1 t \in J_{3l}(u_3)$ for all $t \in[-1,1]$. Then
\begin{eqnarray*}
&&\sum_{l=0}^{L(J)}\int_{J_{3l}(u_3)}
\int_{J_{3l}(u_3)} (K*K)_{h_1}\bigl(x-x'\bigr)
(K*K)_{h_2}\bigl(x-x'\bigr)\\
&&\hspace*{81pt}{}\times f \biggl(x,\frac
{u_3+l}{J}-x
\biggr) \,dx' \,dx
\\
&&\qquad = h_2^{-1}\sum_{l=0}^{L(J)}
\int_{J_{3l}^\mathrm{
o}(u_3;2h_1)}f \biggl(x,\frac{u_3+l}{J}-x \biggr) \,dx \\
&&\hspace*{35pt}\qquad\quad{}\times \int
_{-2}^2 \bigl[K*K(t)\bigr] \bigl[K*K(h_1t/h_2)
\bigr] \,dt +O(1)
\\
& &\qquad= h_2^{-1}\sum_{l=0}^{L(J)}
\int_{J_{3l}(u_3)}f \biggl(x,\frac
{u_3+l}{J}-x \biggr) \,dx\\
&&\hspace*{35pt}\qquad\quad{}\times \int
_{-2}^2 \bigl[K*K(t)\bigr] \bigl[K*K(h_1t/h_2)
\bigr] \,dt +O(1)
\\
&&\qquad = h_2^{-1}f_{w,3}(u_3) \int
_{-2}^2 \bigl[K*K(t)\bigr] \bigl[K*K(h_1t/h_2)
\bigr] \,dt +O(1).
\end{eqnarray*}
This with Lemma~\ref{lem3} below completes the proof of Theorem~\ref{th4}.

\begin{lemma}\label{lem1} Under the condition \textup{(A7)} with the constants $C>0$ and
$\alpha>1/2$,
it follows that \textup{(i)} $J_2^\mathrm{ o}(u_1\dvtx Ch_1^\alpha+h_2) \subset J_2^\mathrm{
o}(u_2;h_2)$
for any $u_1, u_2 \in(a_{k-1}^1,a_k^1)\cap S_1$
with $|u_1-u_2| \le h_1$, $1 \le k \le L_1$; \textup{(ii)} $J_1^\mathrm{
o}(u_1\dvtx Ch_2^\alpha+h_1) \subset
J_1^\mathrm{ o}(u_2;h_1)$
for any $u_1, u_2 \in(a_{k-1}^2,a_k^2)\cap S_2$
with $|u_1-u_2| \le h_2$, $1 \le k \le L_2$.
\end{lemma}

\begin{pf}
We apply (A7) to the choice
$\varepsilon_n =h_1$. Suppose a point $y \in J_2^\mathrm{ o}(u_1; Ch_1^\alpha
+h_2)$. This implies
$y+h_2t +Ch_1^\alpha s \in J_2(u_1)$ for all $s, t \in[-1,1]$. This
holds since $|(h_2t+Ch_1^\alpha s)/(h_2+Ch_1^\alpha)|\le1$
for all $s, t \in[-1,1]$. By (A7), $y+h_2t \in J_2^\mathrm{
o}(u_1;Ch_1^\alpha) \subset J_2(u_2)$ for all $t \in[-1,1]$,
so that we get $y \in J_2^\mathrm{ o}(u_2;h_2)$. The proof of (ii) is the
same.
\end{pf}

\begin{lemma}\label{lem2}
Under the conditions of Theorem~\ref{th4}, It follows that $\bT\hat\mum^A=
o_p(n^{-2/5})$ uniformly on
$S_1 \times S_2 \times S_3$. Furthermore, $\hat\mum^B=n^{-2/5}\tilde
\mum
^B+o(n^{-2/5})$ uniformly
on $S_{1,c}^\mathrm{ o}(h_1)\times S_{2,c}^\mathrm{ o}(h_2)\times S_{3,c}^\mathrm{
o}(C'n^{-\min\{1,\alpha\}/5})$
for a sufficiently large $C'>0$, and $\hat\mum^B(u) =
n^{-2/5}\tilde\mum^B(u) +O(n^{-2/5})$
uniformly on
$S_1 \times S_2 \times S_3$.
\end{lemma}

\begin{pf}
From the standard
theory of
kernel smoothing, it follows that
%
\begin{equation}
\label{bseq} \sup_{(x,y)\in S} \bigl|\hf^A(x,y)\bigr| =
O_p\bigl(n^{-3/10}\sqrt{\log n}\bigr).
\end{equation}
Also, we have $\bA(x,y) =\operatorname{ diag}(1,\nu_2,\nu_2)$ for all $(x,y)$ with
$x \in S_{1,c}^\mathrm{ o}(h_1)$ and $y \in J_2^\mathrm{ o}(x; Ch_1^\alpha+h_2)$,
where $C>0$ and $\alpha>1/2$ are the constants in assumption (A7)
and $\nu_2 = \int u^2 K(u) \,du$.
Define $\cJ=\{(x,y)\in S\dvtx x \in S_{1,c}^\mathrm{ o}(h_1), y \in J_2^\mathrm{
o}(x; Ch_1+h_2)\}$.
From the simplification of $\bA(x,y)$ on $\cJ$, we get
%
\begin{equation}
\label{2Dest} \hf^A(x,y) = \tf^A(x,y),\qquad (x,y) \in\cJ.
\end{equation}
From (\ref{bseq}) and (\ref{2Dest}), we have
%
\begin{eqnarray}
\label{tmu1-1} \hat\mu_1^A(x) = \tmu_1^A(x)+O_p
\bigl(n^{-(3+2r)/10}\sqrt{\log n}\bigr)
\nonumber
\\[-8pt]
\\[-8pt]
\eqntext{\mbox{uniformly for } x \in
S_{1,c}^\mathrm{ o}(h_1),}
\end{eqnarray}
where $r=\min\{1,\alpha\}$. Note that $r>1/2$.
Similarly, we get
%
\begin{eqnarray}
\label{tmu2} \hat\mu_2^A(y) = \tmu_2^A(y)+O_p
\bigl(n^{-(3+2r)/10}\sqrt{\log n}\bigr)
\nonumber
\\[-8pt]
\\[-8pt]
\eqntext{\mbox{uniformly for } y \in
S_{2,c}^\mathrm{ o}(h_2).}
\end{eqnarray}

For the treatment of $\hat\mu_3^A$, we first note that $\bA
(x,(z+l)/J-x)=\operatorname{ diag}(1,\nu_2,\nu_2)$
for all $x \in J_{3l}'(z) \cap S_{1,c}^\mathrm{ o}(h_1)$, where the set
$J_{3l}'(z)$ is defined in the proof of Theorem~\ref{th4}.
In fact,
%
\begin{equation}
\label{seteq} \bigl(x,(z+l)/J-x\bigr) \in\cJ \quad\mbox{if and only if}\quad x \in
J_{3l}'(z) \cap S_{1,c}^\mathrm{ o}(h_1).
\end{equation}
This implies that, for all $0 \le l \le L(J)$,
%
\begin{equation}
\label{2Dest-1}\qquad \hf^A \biggl(x,\frac{z+l}{J}-x \biggr)=
\tf^A \biggl(x,\frac
{z+l}{J}-x \biggr),\qquad x \in
J_{3l}'(z) \cap S_{1,c}^\mathrm{ o}(h_1).
\end{equation}
Due to the condition (A7) we can take a constant $C'>0$ such that,
uniformly for $z \in S_{3,c}^\mathrm{ o}(C'n^{-r/5})$, we have $\sum_{l=0}^{L(J)}\mes[J_{3l}(z)\triangle J_{3l}'(z)]
=O(n^{-r/5})$. Then, from~(\ref{bseq}) and (\ref{2Dest-1}) we have
\begin{eqnarray*}
&&\sum_{l=0}^{L(J)}\int_{J_{3l}(z)}
\hf^A\bigl(x,(z+l)/J-x\bigr) \,dx
\\
&&\qquad = \sum_{l=0}^{L(J)}\int
_{J_{3l}'(z)\cap S_{1,c}^\mathrm{
o}(h_1)}\tf^A\bigl(x,(z+l)/J-x\bigr) \,dx
\\
&&\qquad\quad{}+ O_p\bigl(n^{-3/10}\sqrt{\log n}\bigr) \sum
_{l=0}^{L(J)} \mes\bigl[J_{3l}(z) \triangle
\bigl(J_{3l}'(z)\cap S_{1,c}^\mathrm{ o}(h_1)
\bigr)\bigr]
\\
&&\qquad = \sum_{l=0}^{L(J)}\int
_{J_{3l}(z)}\tf^A\bigl(x,(z+l)/J-x\bigr) \,dx +
o_p\bigl(n^{-2/5}\bigr)
\end{eqnarray*}
uniformly for $z \in S_{3,c}^\mathrm{ o}(C'n^{-r/5})$. This implies
$\hat\mu_3^A(z) = \tmu_3^A(z)+o_p(n^{-2/5})$ uniformly for $z \in
S_{3,c}^\mathrm{ o}(C'n^{-r/5})$.
This together with (\ref{tmu1-1}), (\ref{tmu2}) and Lemma~\ref{lem3} gives
$\bT
\hat\mum^A = o_p(n^{-2/5})$ uniformly on
$S_1 \times S_2 \times S_3$, since $\bT\tilde\mum^A=o_p(n^{-2/5})$
uniformly on the set
and the Lebesgue measures of the set differences $S_1-S_{1,c}^\mathrm{
o}(h_1)$ and $S_2-S_{2,c}^\mathrm{ o}(h_2)$
are of order $n^{-1/5}$ and that of $S_3-S_{3,c}^\mathrm{ o}(C'n^{-r/5})$
is of order $n^{-r/5}$.

To prove the second part of the lemma, recall that $\bA(x,y)=\operatorname{diag}(1,\nu_2,\nu_2)$ on~$\cJ$.
In fact, for $(x,y) \in\cJ$
\[
\int_S \biggl(\frac{u-x}{h_1} \biggr)^j
\biggl(\frac{v-y}{h_2} \biggr)^k K_{h_1}(u-x)K_{h_2}(v-y)
\,du \,dv =0
\]
whenever $j$ or $k$ is an odd integer. This implies $\hf
^B(x,y)=n^{-2/5}\tf^B(x,y) + o(n^{-2/5})$ uniformly for $(x,y)\in\cJ$.
We also get $\hf^B(x,y)=O(n^{-2/5})$ uniformly for $(x,y) \in S$. We
apply the same arguments as in the proof of
the first part, to obtain
\begin{eqnarray*}
\label{tmuB} %
\hat\mu_1^B(x) &=&
n^{-2/5}\tmu_1^B(x) + o\bigl(n^{-2/5}
\bigr) \qquad\mbox{uniformly for } x \in S_{1,c}^\mathrm{ o}(h_1),
\\
\hat\mu_2^B(y) &=& n^{-2/5}\tmu_2^B(y)
+ o\bigl(n^{-2/5}\bigr) \qquad\mbox{uniformly for } y \in S_{2,c}^\mathrm{ o}(h_2).
\end{eqnarray*}
From (\ref{seteq}), it follows that
\begin{eqnarray*}
&&\hf^B \biggl(x,\frac{z+l}{J}-x \biggr)\\
&&\qquad=n^{-2/5}
\tf^B \biggl(x,\frac
{z+l}{J}-x \biggr)+o\bigl(n^{-2/5}
\bigr)
\end{eqnarray*}
for all $(x,z)$ such that $x \in J_{3l}'(z) \cap S_{1,c}^\mathrm{ o}(h_1)$
and $z \in S_3$.
From this and the fact that $\sum_{l=0}^{L(J)}\mes[J_{3l}(z)\triangle
J_{3l}'(z)]=o(1)$
uniformly for $z \in S_{3,c}^\mathrm{ o}(C'n^{-r/5})$, we obtain
\[
\hat\mu_3^B(z) = n^{-2/5}\tmu_3^B(z)
+ o\bigl(n^{-2/5}\bigr) \qquad\mbox{uniformly for } z \in S_{3,c}^\mathrm{ o}
\bigl(C'n^{-r/5}\bigr),
\]
where $C'$ is the constant $C'$ in the proof of the first part.
This completes the proof of the lemma.\vadjust{\goodbreak}
\end{pf}

\begin{lemma}\label{lem3}
Under the conditions of Theorem~\ref{th4}, it follows that
\[
\sup_{u \in S_j}\bigl|\hat\mu_j^A (u)\bigr| =
O_p\bigl(n^{-2/5}\sqrt{\log n}\bigr),\qquad 1 \le j \le3.
\]
\end{lemma}

\begin{pf}
We give the proof for $\hat\mu_1^A$ only. The others are similar.
For $(x,y)$ with $x \in S_1$ and $y \in J_2^\mathrm{ o}(x;Ch_1^\alpha
+h_2)$, we have
\[
\hf^A(x,y) = \varphi_1(x) \hat a_1(x,y) +
\varphi_2(x) \hat a_2(x,y) + \varphi_3(x)\hat
a_3(x,y),
\]
where $\varphi_j$ for $j=1,2,3$ are some bounded functions, $\hat
a_1=\hat b_{00}$, $\hat a_2=\hat b_{10}$ and
$\hat a_3=\hat b_{01}$ with
\begin{eqnarray*}
\hat b_{jk}(x,y) &=& n^{-1}\sum
_{i=1}^n \biggl[ \biggl(\frac
{X_i-x}{h_1}
\biggr)^j \biggl(\frac{Y_i-y}{h_2} \biggr)^k
K_{h_1}(X_i-x)K_{h_2}(Y_i-y)W_i
\\
&&\hspace*{39pt}{} - E \biggl(\frac{X_i-x}{h_1} \biggr)^j \biggl(\frac
{Y_i-y}{h_2}
\biggr)^k K_{h_1}(X_i-x)K_{h_2}(Y_i-y)W_i
\biggr].
\end{eqnarray*}
The lemma follows from (\ref{bseq}) and using
\begin{eqnarray*}
\sup_{x\in S_1} \mes\bigl[J_2(x)-J_2^\mathrm{ o}
\bigl(x;Ch_1^\alpha+h_2\bigr)\bigr] &=&
O_p\bigl(n^{-r/5}\bigr),
\\
\sup_{x \in S_1}\biggl |\int_{J_2(x)} \hat
a_j(x,y) \,dy \biggr| &=& O_p\bigl(n^{-2/5}\sqrt{\log n}
\bigr), \qquad 1 \le j \le3.
\end{eqnarray*}
\upqed\end{pf}
\end{appendix}
%



\printaddresses
\end{document}